\documentclass{article}

\usepackage{arxiv}

\usepackage[utf8]{inputenc} % allow utf-8 input
\usepackage[T1]{fontenc}    % use 8-bit T1 fonts
\usepackage{hyperref}       % hyperlinks
\usepackage{cancel}
\usepackage{url}            % simple URL typesetting
\usepackage{booktabs}       % professional-quality tables
\usepackage{amsfonts}       % blackboard math symbols
\usepackage{nicefrac}       % compact symbols for 1/2, etc.
\usepackage{microtype}      % microtypography
\usepackage{lipsum}		% Can be removed after putting your text content
\usepackage{graphicx}
\usepackage{soul}
\usepackage{doi}
\usepackage{xcolor}
\usepackage{amsmath}
\usepackage{amssymb}
\usepackage{todonotes}
\usepackage{algorithm}
\usepackage{algorithmic}
\graphicspath{{figs/}}
\usepackage{amsthm}
\usepackage{amsfonts}
\usepackage{amsmath}

\newcommand{\prox}{\text{prox}}

\newcommand{\argmin}[1]{\underset{#1}{\text{argmin}}}
\newcommand{\argmax}[1]{\underset{#1}{\text{argmax}}}
\newcommand{\norm}[1]{\lVert #1 \rVert}
\newcommand{\inner}[2]{\langle  #1 , \, #2  \rangle}
\newcommand{\abs}[1]{\vert  #1 \vert}

\title{A Relaxed Primal-Dual Hybrid Gradient Method with Line Search}

\author{ Alex McManus \\
	Department of Applied Mathematics\\
	University of Colorado Boulder\\
	Boulder, CO \\
	\And
        Stephen Becker \\
        Department of Applied Mathematics\\
	University of Colorado Boulder\\
	Boulder, CO \\
        \And
	Nicholas Dwork\\
	Department of Biomedical Informatics\\
	University of Colorado School of Medicine\\
	Aurora, CO
}

\renewcommand{\shorttitle}{Line Search for PDHG}

\begin{document}
\maketitle
\begin{abstract}
The primal-dual hybrid gradient method (PDHG) is useful for optimization problems that commonly appear in image reconstruction. A downside of PDHG is that there are typically three user-set parameters and performance of the algorithm is sensitive to their values. Toward a parameter-free algorithm, we combine two existing line searches. The first, by Malitsky et al. \cite{malitsky2018first}, is over two of the step sizes in the PDHG iterations. We then use the connection between PDHG and the primal-dual form of Douglas-Rachford splitting to construct a line search over the relaxation parameter. We demonstrate the efficacy of the combined line search on multiple problems, including a novel inverse problem in magnetic resonance image reconstruction. The method presented in this manuscript is the first parameter-free variant of PDHG (across all numerical experiments, there were no changes to line search hyperparameters).
\end{abstract}


\section{Introduction}

Convex optimization problems of functions that are not smooth are common in image reconstruction. Proximal algorithms are often used to solve these problems because they are capable of processing optimization variables with many components. This paper is focused on the primal-dual hybrid gradient method (PDHG) by Chambolle and Pock, presented in \cite{chambolle2011first}, which is a first-order proximal algorithm to solve non-smooth problems of the form
\begin{equation*}
%x^{*} \in \underset{x \in \mathcal{H}}{\text{argmin}}\, \left\{f(x) + g(Ax)\right\}.
x^\star \in \underset{x \in \mathcal{H}}{\text{argmin}}\, \left\{f(x) + g(Ax)\right\},
\end{equation*}
where $\mathcal{H}$ is a Hilbert space, $A$ is a linear operator, and $f$ and $g$ are closed convex proper functions.

PDHG can efficiently solve problems that are more general than those that can be solved with other proximal splitting methods (e.g., Douglas-Rachford) because it only requires the proximal operator of $g$, rather than the proximal operator of $g \circ A$ which generally only has a closed form solution when $A A^T$ is a scaled identity.

A downside of PDHG is that, when relaxed, there are three parameters the user needs to set: two step sizes $\tau$ and $\sigma$, and a relaxation parameter $\theta$. Convergence is guaranteed when $\tau\sigma\|A\|^2\leq 1$.  Therefore, if the spectral norm of $A$ is known or can be estimated, the number of parameters that must be chosen by the user can be reduced to two parameters. The rate of convergence is highly dependent on the parameter values which makes a line search that eliminates user-defined parameters desirable. 

Toward that goal, this manuscript presents a combination of two line search techniques over two sets of parameters. Malitsky et al. introduced a line search over the proximal step sizes in \cite{malitsky2018first}. By using the connection between PDHG and the primal-dual formulation of the Douglas-Rachford splitting method (PDDR) \cite{o2020equivalence}, we incorporate the theory of \textit{averaged operator iterations} (AOI).  AOI admits a line search over the relaxation parameter \cite{giselsson2016line}.  The remaining user-set parameter is the ratio between the two scalings for the proximal steps in PDHG.  In general, this can be left as unity, so there are no user-defined parameters remaining.

This paper is organized as follows. Section \ref{sec:preliminaries} contains notation and relevant definitions and facts. Section \ref{sec:algorithms} contains a review of relevant algorithms including the two line search algorithms.  Section \ref{sec:lineSearch} is about the combined line search; its subsection \ref{sec:convergence} details convergence thereof. Section \ref{sec:numExperiments} shows numerical results on several example problems. Section \ref{sec:mri} details an application to magnetic resonance imaging (MRI) reconstruction, a novel combination of two existing MRI reconstruction methods---compressed sensing and partial Fourier reconstruction with homodyne detection. The formulation of the optimization problem for this type of MRI reconstruction is inconvenient to solve with the standard algorithms for compressed sensing reconstruction, which motivated this study of the primal-dual hybrid gradient algorithm.

\section{Preliminaries}
\label{sec:preliminaries}

\subsection{Notation and Definitions}
We follow standard optimization literature notation conventions; see \cite{bauschke2017correction} for further background.
Let $\mathcal{H}$ be a real Hilbert space equipped with inner product $\inner{x}{y}$ and norm $\norm{\,\cdot\,} = \sqrt{\inner{\cdot}{\cdot}}$. Let $\Gamma_0(\mathcal{H})$ be the set of proper, lower semicontinuous convex functions from $\mathcal{H} \to \mathbb{R} \cup \{+\infty\}$. Take $f \in \Gamma_0(\mathcal{H})$. The domain of $f$ is $\text{dom}(f) = \{ x \in \mathcal{H} \, \vert \, f(x) < +\infty\}$. Define as $2^{\mathcal{H}}$ the \textit{power set} of $\mathcal{H}$, or the set of all subsets of the Hilbert space. Let $\mathcal{I}$ denote the identity operator.

The \textit{Legendre-Fenchel conjugate} of function $f$ is $f^* \in \Gamma_0(\mathcal{H})$, given by
\begin{equation*}
  f^*(x) = \sup_{y \in \mathcal{H}}\left\{\inner{x}{y} - f(y)\right\}.
  \label{eq:def_conj}
\end{equation*}
Note that if $f \in \Gamma_0(\mathcal{H})$, then $f^* \in \Gamma_0(\mathcal{H})$ and $f^{**} = f$ \cite{bauschke2017correction}. The \textit{proximal operator} of $f$ is defined as:
\begin{equation*}
  \prox_{\lambda f}(x) = \argmin{u \in \mathcal{H}}\left\{f(u) + \frac{1}{2\lambda}\norm{u - x}^2\right\},
  \label{def:prox}
\end{equation*}
for any $\lambda > 0$.
A \textit{set-valued} operator $A: \mathcal{H} \to 2^{\mathcal{H}}$ assigns to each $x \in \mathcal{H}$ a subset of $\mathcal{H}$ (or an element in $2^{\mathcal{H}}$). The image of $x$ under $A$ is denoted $A(x)$. If $A(x)$ is a singleton, we write $A(x) = y$ instead of $A(x) = \{y\}$. The \textit{graph} of operator $A$ is denoted $\text{gr}(A)$ and defined as
\begin{equation*}
\text{gr}(A) = \left\{(x, y) \, \vert \, y \in A(x)\right\}.
\label{def:graph}
\end{equation*}
The inverse of operator $A$, denoted $A^{-1}$, is defined with graph
\begin{equation*}
\text{gr}(A^{-1}) = \left\{(y, x) \, \vert (x, y) \in \text{gr}(A) \right\}.
\label{def:op_inverse}
\end{equation*}

For operator $A$, define left- and right-scalar multiplication as
\begin{equation*}
\left( \lambda A \right)(x) = \left\{\lambda y \, : \, y \in A(x)\right\}, \qquad \left(A\mu\right)(x) = \left\{y \, : \, y \in A(\mu x)\right\}.
\end{equation*}
The following relationships hold between the graph of $A$ and the graphs of its inverse and scalar multiples:
\begin{equation}
\text{gr}(A^{-1}) = \begin{bmatrix} 0 & 1\\1 & 0\end{bmatrix}\text{gr}(A), \qquad \text{gr}(\lambda A \mu) = \begin{bmatrix} \mu^{-1}\mathcal{I} & 0 \\ 0 & \lambda\mathcal{I}\end{bmatrix}\text{gr}(A).
\label{eq:gr_ops}
\end{equation}

Operator $A$ is \textit{monotone} if
\begin{equation*}\label{def:monotone}
\inner{x - x'}{y - y'} \geq 0, \quad \forall (x, y), (x', y') \in \text{gr}(A).
\end{equation*}
An operator $A$ is \textit{maximally monotone} if it is both monotone and there are no monotone operators $B$ such that $\text{gr}(A) \subset \text{gr}(B)$.

Let $\mathcal{D}$ be a nonempty subset of $\mathcal{H}$ and let $A: \mathcal{D} \to \mathcal{H}$. $A$ is \textit{nonexpansive} if it is Lipschitz continuous with constant 1, or 
\begin{equation*}
\norm{Ax - Ay} \leq \norm{x - y}, \hspace{0.5em} \forall x \in \mathcal{D}, \, \forall y \in \mathcal{D}.
\end{equation*}
$A$ is \textit{firmly nonexpansive} if
\begin{equation*}
\norm{Ax - Ay}^2 + \norm{\left(\mathcal{I} - A\right)x - \left(\mathcal{I} - A\right)y}^2 \leq \norm{x - y}^2, \forall x \in \mathcal{D}, \, \forall y \in \mathcal{D}.
\end{equation*}
We say operator $A$ is \textit{averaged} with constant $\alpha$ if there exists a nonexpansive operator $R: \mathcal{D} \to \mathcal{H}$ such that $A = (1 - \alpha)\mathcal{I} + \alpha R$. $A$ is firmly nonexpansive if and only if it is $1/2-$averaged \cite[Remark 4.34]{bauschke2017correction}.

For a set-valued operator $A$, define its \textit{resolvent} as:
\begin{equation*}
  J_{\lambda A} = (\mathcal{I} + \lambda A)^{-1},
  \label{eq:def_resolv}
\end{equation*}
where $\mathcal{I}$ is the identity. If $A$ is a monotone operator, its resolvent is firmly nonexpansive. Define the \textit{reflected resolvent} of $A$ as:
\begin{equation*}
  R_{\lambda A} = 2J_{\lambda A} - \mathcal{I}.
  \label{eq:def_ref_resolv}
\end{equation*}
Importantly, if $A$ is maximally monotone then its reflected resolvent is nonexpansive \cite[Corollary 23.11]{bauschke2017correction}. The resolvent and reflected resolvents of $A$ satisfy the following identities:
\begin{equation}\label{eq:op_moreau}
    J_{A}(x) + J_{A^{-1}}(x) = x, \quad \text{and} \quad R_{A}(x) + R_{A^{-1}}(x) = 0.
\end{equation}
Since $(\gamma A)^{-1} = A^{-1}\gamma^{-1}$ for scalar $\gamma$ (which can be shown directly from the relationships in \eqref{eq:gr_ops}), 
\begin{equation}\label{eq:op_moreau_scaled}
    J_{\gamma A}(x) + \gamma J_{\gamma^{-1}A^{-1}}(x/\gamma) = x.
\end{equation}
Let $A: \mathcal{H} \to \mathcal{H}$ be such that $AA^* = \alpha\mathcal{I}$, where $\cdot^*$ denotes the operator's adjoint. Let $f \in \Gamma_0(\mathcal{H})$ and let $g(x) = f(Ax)$. Then
\begin{equation*}
\prox_{\gamma g}(x) = x + \frac{1}{\alpha}A^*\left(\prox_{\alpha\gamma f}(Ax) - Ax \right).
\label{def:prox_composed}
\end{equation*}

The \textit{subdifferential} of $f \in \Gamma_0(\mathcal{H})$ is a set-valued operator $\partial f$ defined as:
\begin{equation*}
  \partial f(x) = \left\{\zeta \in \mathcal{H} \, \vert \, f(y) \geq f(x) + \inner{\zeta}{y - x} \, \forall y \in \mathcal{H}\right\}.
  \label{eq:def_subdiff}
\end{equation*}

By the definition of the subdifferential, we have $x^* \in \argmin{x}\, f(x)$ if and only if $0 \in \partial f(x^*)$---this is known as \textit{Fermat's law} \cite[Theorem 16.3]{bauschke2017correction}. Therefore, for suitable functions $f \in \Gamma_0(\mathcal{H})$, solving a minimization problem is the same as solving a \textit{monotone inclusion problem}.

Maximal monotone operators are useful for convex optimization because the subdifferential of $f \in \Gamma_0(\mathcal{H})$ is maximally monotone. Furthermore, the inverse of the subdifferential is the subdifferential of the conjugate function:
\begin{equation*}\label{eq:subdiff_inv}
    (\partial f)^{-1} = \partial f^{*}.
\end{equation*}
An important fact for the convergence of proximal algorithms is that the proximal operator of a proper function $f\in\Gamma_0(\mathcal{H})$ is the resolvent of the subdifferential:
\begin{equation}
\prox_{\gamma f}(x) = \left(\mathcal{I} + \gamma \partial f\right)^{-1}x.
\label{eq:prox_subdiff}
\end{equation}
The subdifferential of $f \in \Gamma_0(\mathcal{H})$ is maximally monotone, so \eqref{eq:prox_subdiff} is single valued.
Letting $A = \partial f$ for $f \in \Gamma_0(\mathcal{H})$ in \eqref{eq:op_moreau_scaled} yields the \emph{Moreau identity}:
\begin{equation*}\label{eq:moreau_prox}    
\prox_{\gamma f}(x) + \gamma\, \prox_{\gamma^{-1}f^*}(x / \gamma) = x.
\end{equation*}

% Another important fact is that for a proper function $f$, i.e. $f \in \Gamma_0(\mathcal{H})$, then the subdifferential of $f$ is a \textit{maximally monotone} operator. One result of this is that the \textit{reflected resolvent} of a maximally monotone operator is nonexpansive which will be important later. \todo{Some of this needs to get spread out or explained better.}

\subsection{The Problem}
For the remainder of this manuscript, let $\mathcal{X}$ and $\mathcal{Y}$ be real Hilbert spaces. We are interested in solving the following (\textit{primal}) optimization problem:
\begin{equation}
%x^* \in \argmin{x \in X}\left\{f(x) + g(Ax) \right\},
x^\star \in \argmin{x \in \mathcal{X}}\hspace{0.5em} f(x) + g(Ax),
\label{eq:prob_def}
\end{equation}
where $f \in \Gamma_0(\mathcal{X})$ and $g \in \Gamma_0(\mathcal{Y})$ have computationally tractable forms of their respective proximal operators and $A: \mathcal{X} \to \mathcal{Y}$ is a linear operator with known adjoint $A^*$ and finite operator norm:
\begin{equation*}
  \norm{A} = \sup\left\{ \norm{Ax} \, : \, x \in \mathcal{X} \, \text{with } \norm{x} \leq 1 \right\} < +\infty.
\label{eq:a_norm}
\end{equation*}
% \begin{itemize}
%     \item $f \in \Gamma_0(\mathcal{X})$ and $g \in \Gamma_0(\mathcal{Y})$ have closed forms of their respective proximal operators, and
    % \item $A: \mathcal{X} \to \mathcal{Y}$ is a bounded linear operator with known adjoint $A^*$ and finite operator norm
    % \begin{equation*}
    %     \norm{A} = \sup\left\{ \norm{Ax} \, : \, x \in \mathcal{X} \, \text{with } \norm{x} \leq 1 \right\} < +\infty.
    % \label{eq:a_norm}    
    % \end{equation*}
% \end{itemize} 
The corresponding dual problem to \eqref{eq:prob_def} is:
\begin{equation}
%y^* \in \argmax{y \in \mathcal{Y}} -\left\{f^*(-A^*y) + g^*(y)\right\}.
y^\star \in \argmax{y \in \mathcal{Y}}\hspace{0.5em} ^-\left(f^*(-A^*y) + g^*(y)\right).
\label{eq:dual_prob}
\end{equation}
The algorithm that is the focus of this manuscript is \textit{primal-dual}, meaning it is applied to the saddle point interpretation of the primal \eqref{eq:prob_def} and dual \eqref{eq:dual_prob} problems:

\begin{equation}
(\hat{x}, \hat{y}) \in \argmin{x \in \mathcal{X}}\;\argmax{y \in \mathcal{Y}} \left\{\inner{Ax}{y} + f(x) - g^*(y)\right\}.
\label{eq:saddle}
\end{equation}
More details can be found in, e.g., \cite{chambolle2011first, condat2013primal, bauschke2017correction}.

The saddle point formulation yields the following optimality conditions, which are the starting point for development of primal-dual algorithms to solve the problem:

\begin{equation}
(\hat{x}, \hat{y}) \in \mathcal{X} \times \mathcal {Y} \;\; \text{such that } \begin{pmatrix} 0 \\ 0 \end{pmatrix} \in \begin{pmatrix} 0 & A^* \\ -A & 0\end{pmatrix} \begin{pmatrix}x \\ y \end{pmatrix} + \begin{pmatrix} \partial f(x) \\ \partial g^*(y)\end{pmatrix}.
\label{eq:opt_cond}
\end{equation}

The Karush-Kuhn-Tucker conditions state that if $(\hat{x}, \hat{y})$ is a solution to \eqref{eq:opt_cond}, then $\hat{x}$ is a solution to \eqref{eq:prob_def} and $\hat{y}$ is a solution to \eqref{eq:dual_prob} \cite[Theorem 19.1]{bauschke2017correction}. Furthermore, $(\hat{x}, \hat{y})$ is a solution to \eqref{eq:saddle} \cite[Proposition 19.20]{bauschke2017correction}. Note that the converse is not true in general and the set of solutions to \eqref{eq:opt_cond} may be empty. Under certain constraint qualifications such as
\begin{equation}
0 \in \text{sri}\left(\text{dom}(g) - A(\text{dom}(f))\right),
\label{eq:constraint_qualification}
\end{equation}
then the duality gap is zero (i.e. strong duality holds) and the set of solutions to \eqref{eq:dual_prob} is nonempty \cite[Theorem 15.23]{bauschke2017correction}. Here $\text{sri}(C)$ is the \textit{strong relative interior} of a set $C$ (see \cite[Section 6.2]{bauschke2017correction} for more details). In many applications of interest, $\mathcal{X}$ and $\mathcal{Y}$ are finite-dimensional. The condition \eqref{eq:constraint_qualification} is then satisfied if
\begin{equation}
(\text{ri dom}(g)) \cap A(\text{ri dom}(f)) \neq \varnothing,
\label{eq:constraint_qualification2}
\end{equation}
where $\text{ri}$ is the \textit{relative interior} of a set \cite[Proposition 15.24]{bauschke2017correction}. The condition \eqref{eq:constraint_qualification2} is easy to show if either $f$ or $g$ is full domain. Under this condition, whenever $\hat{x}$ is a solution to \eqref{eq:prob_def} and $\hat{y}$ is a solution to \eqref{eq:dual_prob}, then $(\hat{x}, \hat{y})$ is a solution to \eqref{eq:opt_cond}.

\section{Algorithms}
\label{sec:algorithms}
\subsection{Douglas-Rachford Splitting}
The standard form of Douglas-Rachford splitting (DR) was originally given by Lions and Mercier \cite{lions1979splitting} to find zeros of a sum of maximally monotone operators $F$ and $G$; i.e., to find $x$ such that
\begin{equation*}
0 \in F(x) + G(x).
\label{eq:dr_monotone}
\end{equation*}
The standard version of the algorithm is the fixed point iteration
\begin{equation}
y_{k+1} = \frac{1}{2}y_k + \frac{1}{2}R_{\gamma G}R_{\gamma F}y_k,
\label{eq:drWithoutRelax}
\end{equation}
where $\gamma$ is a positive constant.

This is commonly presented as the following steps:
\begin{equation*}
\begin{aligned}
x_{k+1} &= J_{\gamma F}(y_k) \\
w_{k+1} &= J_{\gamma G}(2x_{k+1} - y_k)\\
y_{k+1} &= y_k + w_{k+1} - x_{k+1}.
\end{aligned}
\end{equation*}
If we allow for a relaxation parameter $\alpha_k \in (0, 1)$, then \eqref{eq:drWithoutRelax} is altered as follows:
\begin{equation}
y_{k+1} = (1 - \alpha_k)y_k + \alpha_kR_{\gamma G}R_{\gamma F}y_k.
\label{dr_aoi}
\end{equation}
Equivalently,
\begin{equation}
\begin{aligned}
x_{k+1} &= J_{\gamma F}(y_k) \\
w_{k+1} &= J_{\gamma G}(2x_{k+1} - y_k)\\
y_{k+1} &= y_k + \rho_k(w_{k+1} - x_{k+1}),
\end{aligned}
\label{eq:dr_relax_iters}
\end{equation}
where $\rho_k = 2\alpha_k$.

Results on the convergence of DR are given in, e.g., \cite{bauschke2017correction, lions1979splitting, bauschke2016douglas, eckstein1992douglas}. The standard assumptions for convergence are the maximal monotonicity of $F$ and $G$ and the existence of a solution. Using relationship \eqref{eq:op_moreau} on the second step of the iteration \eqref{eq:dr_relax_iters} yields the iterations for \textit{primal dual} Douglas-Rachford (PDDR):

\begin{equation}
\begin{aligned}
    x_{k+1} &= J_{\gamma F}(y_k) \\
    w_{k+1} &= J_{(\gamma G)^{-1}}(2x_{k+1} - y_k) \\
    y_{k+1} &= (1 - \rho_k)y_k + \rho_k(x_{k+1} - w_{k+1}).
\end{aligned}
\label{eq:iterations_pddr_resolv}
\end{equation}

By using relationship \eqref{eq:op_moreau} applied to the fixed point iteration version of DR in \eqref{dr_aoi}, we get the fixed point iteration form of PDDR:
\begin{equation}
y_{k+1} = (1 - \alpha_k)y_k -\alpha_k R_{(\gamma G)^{-1}}R_{\gamma F}y_k.
\label{eq:pddr_aoi}
\end{equation}

\subsection{Primal-Dual Hybrid Gradient}

The Primal-Dual Hybrid Gradient (PDHG) method that solves problems of the form of \eqref{eq:opt_cond} was originally proposed by Pock et al. in \cite{pock2009algorithm}. Convergence of the algorithm was formalized by Chambolle and Pock in \cite{chambolle2011first}, so it is alternatively called the \textit{primal-dual hybrid gradient} method (PDHG) or the \textit{Chambolle-Pock} algorithm. The method has the following iterations:
\begin{equation}
\begin{aligned}
\bar{x}_{k+1} &= J_{\tau_k F}\left(x_k - \tau_k A^* z_k\right) \\
\bar{z}_{k+1} &= J_{\sigma_k G^{-1}}\left(z_k + \sigma_k A(2\bar{x}_{k+1} - x_k)\right) \\
x_{k+1}&=x_{k} + \alpha_k(\bar{x}_{k+1} - x_k) \\
z_{k+1}&=z_{k} + \alpha_k(\bar{z}_{k+1} - z_k),
\end{aligned}
\label{eq:iterations-pdhg}
\end{equation}
where $F = \partial f$, $G^{-1} = \partial g^*$ for proper, lower semicontinuous convex functions $f \in \Gamma_0(\mathcal{X})$ and $g \in \Gamma_0(\mathcal{Y})$, and $A : \mathcal{X} \to \mathcal{Y}$ is a bounded linear operator with known adjoint. Since $g \in \Gamma_0(\mathcal{Y})$, $g^* \in \Gamma_0(\mathcal{Y})$ as well. The resolvents $J_F$ and $J_{G^{-1}}$ are proximal operators of the functions $f$ and $g$, as in \eqref{eq:prox_subdiff}.

As written, the algorithm has three parameters for each iteration: $\tau$ and $\sigma$, which must satisfy $\tau\sigma\norm{A}^2 \leq 1$ to guarantee convergence \cite{condat2013primal}, and $\alpha$, which is a \textit{relaxation parameter} and typically takes values in $(0, 2)$ \cite{bauschke2017correction}.  Commonly, parameter values are held constant for all iterations.  I.e., $\tau_k=\tau$, $\sigma_k=\sigma$, and $\alpha_k=\alpha$ for all $k$.

A line search for PDHG was introduced in \cite{malitsky2018first} by Malitsky et al. The pseudocode for the line search is show in Alg. \ref{alg:pdhgls}. Line searches for primal-dual methods are difficult in general. The conditions for a standard line search, e.g., Armijo backtracking line search \cite{armijo1966minimization}, cannot always be evaluated (e.g., when primal-dual methods are applied to problems where $f$ or $g$ are an indicator function).

\begin{algorithm}[!t]\caption{\textit{PDHG Line Search }\cite{malitsky2018first}}
    \label{alg:pdhgls}
\begin{algorithmic}
   \STATE {\bfseries Initialization:}
Given $x_{k-1}, z_{k-1} \in \mathcal{H}$, $\prox_f$, $\prox_{g^*}$, $\tau_{k-1} > 0$, $\theta_{k - 1} > 0$, $\mu \in (0, 1)$, $\beta > 0$, $\delta \in (0,1)$.
\STATE {\bfseries Main iteration:}
\STATE 1. Compute  
\begin{align*}
    x_k= \prox_{\tau_{k-1}f}(x_{k-1} - \tau_{k-1}A^{*}z_{k-1})
\end{align*}
\STATE 2. Choose $\tau_{k} \in [\tau_{k-1}, \tau_{k-1}\sqrt{1+\theta_{k-1}}]$ and run \hfill\COMMENT{By default we choose $\tau_k = \tau_{k-1}\sqrt{1+\theta_{k-1}}$}\\
\STATE ~~~ \textbf{Linesearch:}
\STATE ~~~ 2.a. Compute
\vspace{-3ex}
\begin{align*}
  \theta_{k} & = \frac{\tau_{k}}{\tau_{k-1}}\\
  \bar{x}_k &= x_k + \theta_k\left(x_k - x_{k-1}\right)\\
  z_{k} & = \prox_{\beta \tau_{k} g^*}(z_{k-1} + \beta\tau_{k} A\bar{x}_{k})
\end{align*}
    \STATE   ~~~2.b. Break linesearch if
\begin{equation*}\label{stop_crit_s}
    \sqrt{\beta}\tau_{k+1} \lVert A^*z_k - A^*z_{k-1}\rVert \leq \delta \lVert z_k - z_{k-1}\rVert
\end{equation*}
and return
\begin{equation*}
z_k, \, \bar{x}_k
\end{equation*}
\vspace{-3ex}
\STATE ~~~~~ \quad Otherwise, set $\tau_k\gets \mu\tau_k $ and go to 2.a.
\STATE ~~~\textbf{End of linesearch}

\end{algorithmic}
\end{algorithm}

\subsection{The Relationship Between PDDR and PDHG}

It has been shown in \cite{o2020equivalence} that the PDDR iterations \eqref{eq:iterations_pddr_resolv} applied to the following \textit{modified problem} of \eqref{eq:DRtoPDHG} yield the same iterations as PDHG applied to the original problem \eqref{eq:prob_def}:
\begin{equation}
    (\hat{x}, \hat{y}) \in \argmin{x \in \mathcal{X}, y \in \mathcal{Y}} \tilde{f}(x, y) + \tilde{g}(x, y)
\label{eq:DRtoPDHG}
\end{equation}
with modified functions
\begin{equation*}
\tilde{f}(x, y) = f(x) + \delta_{\{0\}}(y), \qquad \tilde{g}(x, y) = g(Ax + By),
\end{equation*}
where $\delta_{\{0\}}$ is an indicator function specifying that the argument is $0$ and $B$ is chosen to satisfy $AA^* + BB^* = \frac{1}{\theta}\mathcal{I}$ for some $\theta > 0$.  (Note that an indicator takes on a value of $0$ when its argument satisfies its condition; otherwise, it takes on a value of $\infty$.) As an example, if $A$ is an $n \times n$ matrix with $\frac{1}{\theta} \geq \norm{A}_2^2$, let $B = \left(\frac{1}{\theta}\mathcal{I} - AA^*\right)^{1/2}$ where $\left(\cdot\right)^{1/2}$ is the Cholesky factorization. Defining these modified functions and applying the iterations of PDDR \eqref{eq:iterations_pddr_resolv} to problem \eqref{eq:DRtoPDHG} yields the PDHG iterations of \eqref{eq:iterations-pdhg} \cite{o2020equivalence}. In \cite{o2020equivalence}, the choice of $\theta$ is left to the user. The relationship between $\theta$ and the $\tau, \sigma$ parameters from PDHG is $\tau\sigma = \theta$ \cite{o2020equivalence}. 

As can be seen in \cite{o2020equivalence}, the relationship between the variables of PDDR in \eqref{eq:iterations_pddr_resolv} when applied to \eqref{eq:DRtoPDHG} and PDHG in \eqref{eq:iterations-pdhg} when applied to \eqref{eq:prob_def} is
\begin{align}
y_{k} = \left(x_k - \tau_k A^*z_k, -\tau_k B^*z_k\right).
\label{eq:pddr_to_pdhg}
\end{align}
This relationship is what allows us to tie the two existing line searches together. An iteration of PDHG is the same as one iteration of PDDR with the relationship above.

\subsection{The Averaged Operator Iteration Line Search}
An \textit{averaged operator iteration} (AOI) for a nonexpansive operator $S$ is an iteration of the form
\begin{equation}
y_{k+1} = (1 - \alpha)y_k + \alpha S y_k = y_k + \alpha(Sy_k - y_k).
\label{eq:def_aoi}
\end{equation}
All of the algorithms discussed in this section so far can be written as an AOI. The AOI form of Douglas-Rachford is given in \eqref{eq:drWithoutRelax} where $S = R_{\gamma G}R_{\gamma F}$. The AOI form of PDDR is given in \eqref{eq:pddr_aoi} with $S = -R_{(\gamma G)^{-1}}R_{\gamma F}$. In both cases, $S$ is nonexpansive as it is a composition of reflected resolvents, which are themselves nonexpansive \cite[Proposition 4.4]{bauschke2017correction}.

A line search for averaged operator iterations was introduced by Gisselson et al.\ in \cite{giselsson2016line}. The pseudocode for the algorithm is presented in Alg. \ref{alg:aoi}. 

\begin{algorithm}[!t]\caption{\textit{AOI Line Search \cite{giselsson2016line}}}
    \label{alg:aoi}
\begin{algorithmic}
   \STATE {\bfseries Initialization:}
Choose $y_0\in \mathcal{H}$, $\alpha_{\text{max}} > 0$, $\bar{\alpha} \in (0, \alpha_{\text{max}})$, $\varepsilon > 0$, $\mu \in (0,1)$.
\STATE {\bfseries Main iteration:}
\STATE 1. Compute  
\begin{align*}
    r_{k} &= Sy_{k} - y_{k}\\
\bar{y}_k &= y_k + \bar{\alpha}r_k\\
\bar{r}_k &= S\bar{y}_{k} - \bar{y}_{k}
\end{align*}
\STATE 2. Set $\alpha_k = \alpha_{\text{max}}$ and run\\
\STATE ~~~ \textbf{Linesearch:}
\STATE ~~~ 2.a. Compute
\vspace{-3ex}
\begin{align*}
  y_{k+1} &= y_k + \alpha_kr_k\\
  r_{k+1} & = Sy_{k+1} - y_{k+1}\\
\end{align*}
	\vspace{-4ex}
    \STATE   ~~~2.b. Break linesearch if
\begin{equation*}\label{stop_crit}
    \lVert{r_{k+1}\rVert} = \lVert{Sy_{k+1} - y_{k+1}\rVert} \leq (1-\varepsilon)\lVert{\bar{r}_k\rVert}
\end{equation*}
\vspace{-3ex}
\STATE ~~~~~ \quad Otherwise, set $\alpha_k\gets \mu\alpha_k $ and go to 2.a. If $\alpha_k = \bar{\alpha}$, terminate the line search.
\STATE ~~~\textbf{End of linesearch}
\end{algorithmic}
\end{algorithm}

AOI with line search requires the user to set a nominal step size $\bar{\alpha}$. This can be set to $\frac{1}{2}$ to ensure the \textit{nominal} step corresponds to a fixed-point iteration of a firmly nonexpansive operator. In order, the steps of Alg. \ref{alg:aoi} compute the current residual, compute the \textit{nominal} next step with the nominal step size, and then the residual with the nominal step size. The line search then iterates over the values of $\alpha_k$ in the final step. Either the step size is $\alpha_k = \bar{\alpha}$, the nominal step, or $\alpha_k \in (\bar{\alpha}, \alpha_{\text{max}}]$ is chosen such that
\begin{equation}
    \norm{r_{k+1}}_2 = \norm{Sy_{k+1} - y_{k+1}}_2 \leq (1 - \varepsilon)\norm{\bar{r}_{k}}_2,
    \label{eq:aoi_linesearch_cond}
\end{equation}
where $\varepsilon > 0$ and $\alpha_{\text{max}} > \bar{\alpha}$ are fixed algorithm parameters. 

The AOI iterations of \eqref{eq:def_aoi} can be interpreted as steps of length $\alpha$ in the direction of the residual $Sx_k - x_k$. Taking a large step in this direction may be beneficial and that is what this line search aims to accomplish. In practice, $\alpha_{\text{max}}$ is set much larger than $\bar{\alpha}$ and progressively shrunk to test if a large step in the residual direction reduces the norm of the residual more than the nominal step.

\section{The Line Search}
\label{sec:lineSearch}
The connection between PDDR as an AOI has been shown and the relationship between the variables of PDDR and PDHG are presented in \eqref{eq:pddr_to_pdhg}. Thus the strategy of the proposed line search:
\begin{itemize}
    \item Use the AOI line search to find an advantageous relaxation parameter, where for each proposed relaxation parameter we
    \item convert the variables to those of PDHG and use the line search of Malitsky et al. to find proximal step sizes $\tau$ and $\sigma$.
\end{itemize}
We refer to these steps, respectively, as the ``outer'' and ``inner'' line searches, as the AOI line search calls the PDHG line search at each step. 
We propose a \textit{relaxed primal-dual hybrid gradient method with line search}, or rPDHG. The line search executes the Malitsky line search for PDHG \cite{malitsky2018first} and uses the relationship between PDHG and the primal-dual form of Douglas-Rachford to create a line search over the relaxation parameter.  
Pseudocode for the Outer Line Search is shown in Alg. \ref{alg:A1}, and pseudocode for the Inner Line Search (which is equivalent to the Maltisky PDHG line search) is shown in Alg. \ref{alg:A2}.

The proposed line search inherits some parameters from the AOI and PDHG line searches. The $\beta$ parameter is the ratio between the two proximal step sizes $\tau_k$ and $\sigma_k$, which we leave as unity in experiments but include in the algorithm for completeness.

\begin{algorithm}[!t]\caption{\textit{Relaxed PDHG with Line Search - Outer Line Search}}
    \label{alg:A1}
\begin{algorithmic}
   \STATE {\bfseries Initialization:}
Choose $x_0, z_0 \in \mathcal{H}$, $\alpha_{\text{max}} > 0$, $\bar{\alpha} \in (0, \alpha_{\text{max}})$, $\varepsilon > 0$, $\mu \in (0,1)$, $\tau_0 > 0$, $A, B$ such that $AA^*+BB^* = \frac{1}{\theta}\mathcal{I}$.
\STATE {\bfseries Main iteration:}
\STATE 1. Compute
\begin{equation*}
[x_{k+1/2}, z_{k+1/2}, \tau_{k}] = \text{Inner Line Search}(x_{k-1}, z_{k-1}, \tau_{k-1}, \bar{\alpha}).
\end{equation*}
\STATE 2. Compute
\begin{align*}
\bar{x}_{k+1} &= \prox_{\tau_{k} f}\left(x_{k+1/2}-\tau_k A^*z_{k+1/2}\right),\\
\bar{z}_{k+1} &= \prox_{\beta\tau_{k}g^*}\left(z_{k+1/2} + \beta\tau_{k}A(2*\bar{x}_{k+1} - x_{k+1/2})\right),
\end{align*}
\STATE 3. Compute
\begin{equation*}
\bar{r}_k = \begin{bmatrix} \left(\bar{x}_k - \tau_kA^*\bar{z}_{k}\right) - \left(x_{k+1/2} - \tau_{k+1/2}A^*z_{k+1/2}\right) \\ -\tau_{k}B^*\bar{z}_k - \left(-\tau_{k+1/2}B^*z_{k+1/2}\right) \end{bmatrix}.
\end{equation*}
\STATE 4. Set $\alpha_k = \alpha_{\text{max}}$ and run\\
\STATE ~~~ \textbf{Linesearch:}
\STATE ~~~ 4.a. Compute
\vspace{-2ex}
\begin{align}
[x_{k+1/2}, z_{k+1/2}, \tau_{k}] &= \text{Inner Line Search}(x_{k-1}, z_{k-1}, \tau_{k-1}, \alpha_k),\\
\bar{x}_{k+1} &= \prox_{\tau_{k} f}\left(x_{k+1/2}-\tau_k A^*z_{k+1/2}\right),\\
\bar{z}_{k+1} &= \prox_{\beta\tau_{k}g^*}\left(z_{k+1/2} + \beta\tau_{k}A(2*\bar{x}_{k+1} - x_{k+1/2})\right),\\
r_{k+1} &= \begin{bmatrix} \left(\bar{x}_{k+1} - \tau_kA^*\bar{z}_{k+1}\right) - \left(x_{k+1/2} - \tau_{k}A^*z_{k+1/2}\right) \\ -\tau_{k}B^*\bar{z}_{k+1} - \left(-\tau_{k}B^*z_{k+1/2}\right) \end{bmatrix}.
\end{align}
	\vspace{-2ex}
    \STATE   ~~~4.b. Break line search if
\begin{equation*}\label{gpdhg:stop_crit}
    \lVert{r_{k+1}\rVert} \leq (1-\varepsilon)\lVert{\bar{r}_k\rVert}
\end{equation*}
\vspace{-3ex}
\STATE ~~~~~ \quad Otherwise, set $\alpha_k\gets \mu\alpha_k $ and go to 4.a. If $\alpha_k = \bar{\alpha}$, terminate the line search.
\STATE ~~~\textbf{End of line search}
\end{algorithmic}
\end{algorithm}

\begin{algorithm}[!ht]\caption{\textit{Relaxed PDHG with Line Search - Inner Line Search}}
    \label{alg:A2}
\begin{algorithmic}
   \STATE {\bfseries Initialization:}
Given $x_{k-1}, z_{k-1}\in \mathcal{H}$, $\prox_f$, $\prox_{g^*}$, $\alpha_k \in [\bar{\alpha}, \alpha_{max}]$, $\tau_{k-1} > 0$, $\theta_{k - 1} > 0$, $A$, $B$ such that $AA^* + BB^* = \frac{1}{\theta}\mathcal{I}$, $\mu \in (0, 1)$, $\beta > 0$, $\delta \in (0,1)$.
\STATE {\bfseries Main iteration:}
\STATE 1. Compute  
\begin{align*}
   \hat{x}_k= \prox_{\tau_{k-1}f}(x - \tau_{k-1}A^{*}z_{k-1})
\end{align*}
\STATE 2. Choose $\tau_{k} \in [\tau_{k-1}, \tau_{k-1}\sqrt{1+\theta_{k-1}}]$ and run \hfill\COMMENT{By default we choose $\tau_k = \tau_{k-1}\sqrt{1+\theta_{k-1}}$}\\
\STATE ~~~ \textbf{Line search:}
\STATE ~~~ 2.a. Compute
\vspace{-3ex}
\begin{align*}
  \theta_{k} & = \frac{\tau_{k}}{\tau_{k-1}}\\
  \bar{x}_k &= \hat{x}_k + \theta_k\left(\hat{x}_k - x_{k-1}\right)\\
  \bar{z}_{k} & = \prox_{\beta \tau_{k} g^*}(z_{k-1} + \beta\tau_{k} K\bar{x}_{k})
\end{align*}
    \STATE   ~~~2.b. Break line search if
\begin{equation*}\label{gpdhg:stop_crit_s}
    \sqrt{\beta}\tau_{k} \lVert A^*\bar{z}_k - A^*z_{k-1}\rVert \leq \delta \lVert \bar{z}_k - z_{k-1}\rVert.
\end{equation*}
\STATE ~~~~~ \quad Otherwise, set $\tau_k\gets \mu\tau_k $ and go to 2.a.
\STATE ~~~\textbf{End of line search}
\STATE 3. Compute
\begin{align*}
x_k &= (1 - 2\alpha_k)x_{k-1} + 2\alpha_k\bar{x}_k\\
z_k &= (1 - 2\alpha_k)z_{k-1} + 2\alpha_k\bar{z}_{k}.
\end{align*}
\STATE Return $x_k, z_k, \tau_k, \theta_k$. 
%\vspace{-3ex}
\end{algorithmic}
\end{algorithm}

\subsection{Convergence}
\label{sec:convergence}
Convergence of the outer (AOI) line search is detailed in the original paper \cite{giselsson2016line}. To briefly summarize, let $k_0$ denote the number of times the outer line search activates (that is, chooses $\alpha_k > \bar{\alpha}$). If $k_0$ is finite, let $k_{\text{max}}$ denote the final iteration on which $\alpha_k > \bar{\alpha}$. For all $k > k_{\text{max}}$, the iteration is 
\begin{equation*}
y_{k+1} = y_{k} + \bar{\alpha}\left(S_ky_k - y_k\right).
\label{eq:conv1}
\end{equation*}
For $\bar{\alpha} \in (0, 1)$ this corresponds exactly to an iteration of PDHG with Malitsky's line search. As shown in \cite[Lemma 1]{malitsky2018first}, this line search always terminates. Convergence is equivalent to convergence of PDHG with line search.

For each iteration $k$ that $\alpha_k > \bar{\alpha}$ the line search condition requires that (for $\varepsilon \in (0, 1)$)
\begin{equation*}
\norm{r_{k+1}} = \norm{S_ky_k - y_k} \leq (1 - \varepsilon)\norm{\bar{r}_k} \leq (1 - \varepsilon)\norm{r_k},
\label{eq:conv2}
\end{equation*}
which yields
\begin{equation*}
\norm{r_{k+1}} \leq (1 - \varepsilon)^{k_0}\norm{r_0}.
\label{eq:conv3}
\end{equation*}
If $k_0 \to \infty$ as $k \to \infty$, then $\norm{r_{k+1}} \to 0$.

\subsection{Line Search Activation}\label{sec:ls_activation}
Conducting the outer line search is computationally costly. If there were an inexpensive test that could be run to determine whether or not line search would be beneficial, then one could activate line search only in these instances to save time.  In \cite{giselsson2016line}, Gisselson et al. suggest a test based on the displacement vectors of the trajectory: $v_k = y_k - y_{k-1}$.  They suggest activating the line search when $v_k$ is approximately aligned with $v_{k+1}$; i.e., when the following is true
\begin{equation*}
  \frac{ v_{k}^T v_{k-1} }{ \| v_{k} \|_2 \| v_{k-1} \|_2 } > 1 - \hat{\varepsilon}
\end{equation*}
for some small $\hat{\varepsilon}>0$.

We hypothesize that this heuristic was useful because the objective function of the problem in their experiments was a strongly-convex quadratic where the gradient becomes smaller quickly as one approaches the optimal point.  The $\ell_1$ norm used as the objective function in our application does not have this property, and the line search was almost always activated.

Instead, we propose the following: the line search will be activated if either 1) it is the first iteration, 2) line search was activated in the previous iteration and $\alpha_k$ did not equal $\bar{\alpha}$, or 3) $\|r_{k}\|_2 / \| r_{k-1} \|_2 < 1 - \hat{\varepsilon}$ for some small $\hat{\varepsilon}>0$.  For the results presented in this manuscript, $\hat{\varepsilon}=0.05$.  In words, the third condition means that line search should be activated when the size of the residual for the current iteration is smaller than that of the previous iteration.  We based this on \eqref{eq:aoi_linesearch_cond}, where we are using $r_{k-1}$ as a surrogate for $\bar{r}_k$.

\subsection{Drawbacks}\label{sec:drawbacks}
The line search presented may not be suitable for every application. For example, in image processing applications $A$ may operate on $n \times n$ images, so $A: \mathbb{R}^{n^2} \to \mathbb{R}^{n^2}$. For modest $n$ this can quickly exceed memory limitations for standard computers if $A$ is stored in dense matrix form. Many of these operators work matrix-free, but finding a suitable $B$ for the proposed line search may require storing $B$ as a dense matrix where $B$ is at least as large as $A$. An example of this phenomenon is shown in Section \ref{sec:2dtv}, where the image had to be scaled down so $B$ could fit in memory.

When $A$ is matrix free, it may be possible to find a matrix-free implementation of $B$; an example is shown in Section \ref{sec:mri}. Alternatively, one may be able to evaluate $B$ matrix-free for arbitrary $A$ by using the Lanczos method \cite{amsel2024}, but it is unclear whether this would work in practice due to the tradeoff between accuracy and cost.
\section{Numerical Experiments}
\label{sec:numExperiments}
We present numerical results for the following problems:
\begin{itemize}
	\item A generalized LASSO-style regularized least squares problem,
	\item a one-dimensional (signal) total variation denoising problem,
	\item a two-dimensional (image) total variation denoising problem, and
	\item a novel magnetic resonance image (MRI) reconstruction problem using compressed sensing.
\end{itemize}

For the numerical experiments, we compare the proposed line search against:
\begin{itemize}
	\item Standard PDHG iterations \eqref{eq:iterations-pdhg},
	\item PDHG with the line search of Malitsky et al. as in Alg. \ref{alg:pdhgls},
	\item PDDR applied to the modified problem using the AOI formulation and line search as in Alg. \ref{alg:aoi}.
\end{itemize}
Extensive parameter searches were conducted to choose the optimal parameters for these algorithms, used in the comparison to the proposed method.

\subsection{LASSO}
\begin{figure}[t]
\centering
\includegraphics[width=0.75\linewidth]{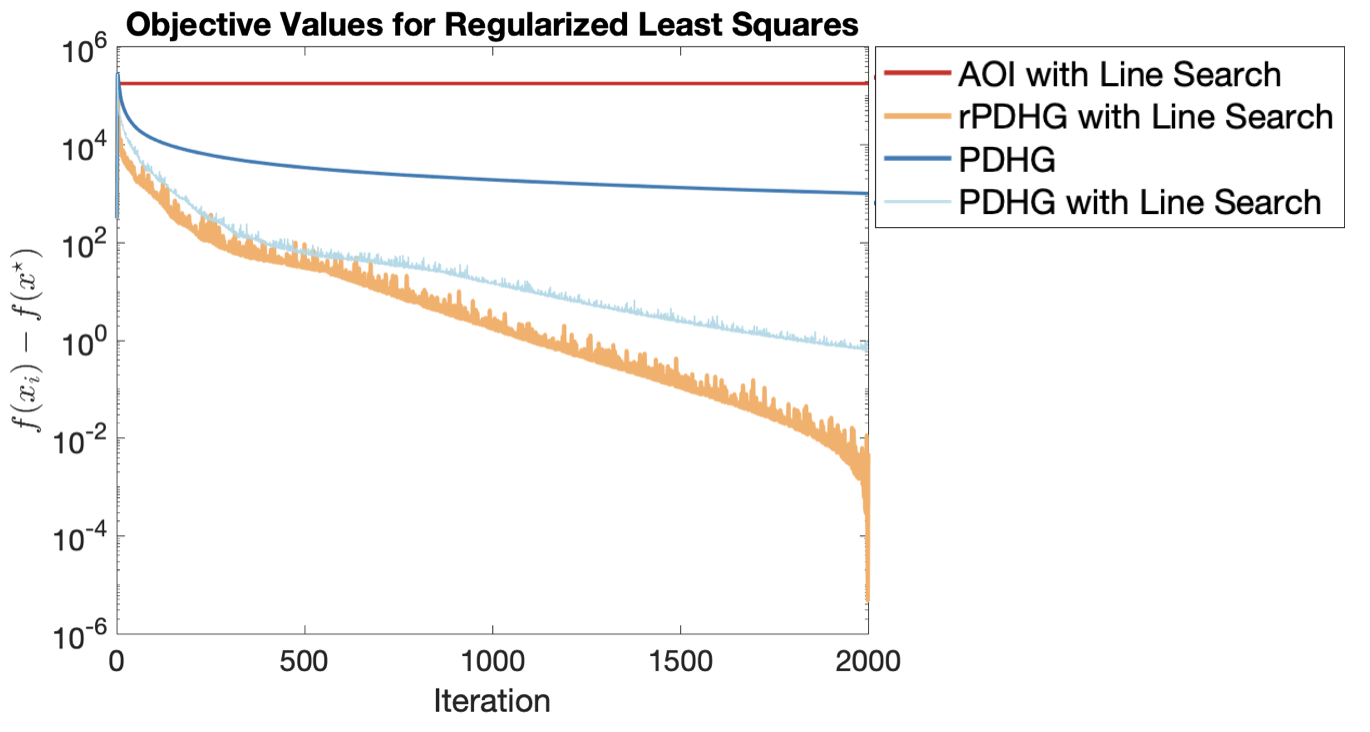}
\caption{A comparison of the objective values from the different optimization algorithms as they solve the regularized least squares problem in \eqref{res:lasso}. }
\label{fig:lasso}
\end{figure}

Consider the unconstrained generalized LASSO problem:
\begin{equation}
\underset{x \in \mathbb{R}^n}{\text{minimize}}\, \frac{1}{2}\lVert x - b \rVert_2^2 + \lambda \lVert Ax \rVert_1.
\label{res:lasso}
\end{equation}
Here $A \in \mathbb{R}^{n\times n}$, $x \in \mathbb{R}^n$, and $b \in \mathbb{R}^n$. We write $f(x) = \frac{1}{2}\lVert x - b \rVert_2^2$ and $g(x) = \lambda \lVert x \rVert_1$. $A$ and $b$ are randomly generated from a standard normal distribution. $B$ (for the relaxed PDHG with line search) is chosen by setting $\theta = 0.9 / \norm{A}_2^2$ and solving for $B$ from $BB^* = \frac{1}{\theta}\mathcal{I} - AA^*$ via Cholesky factorization. For the results presented, we set $n = 1000$. After searching over parameters, we compare the best convergence rate for each method. This is shown in Fig. \ref{fig:lasso}. For rPDHG, however, we set $\beta = 1$ and choose a $\tau_0 > 0$, avoiding a costly parameter search.

\subsection{1D-TV Denoising}

\begin{figure}[t]
\centering
\includegraphics[width=0.5\linewidth]{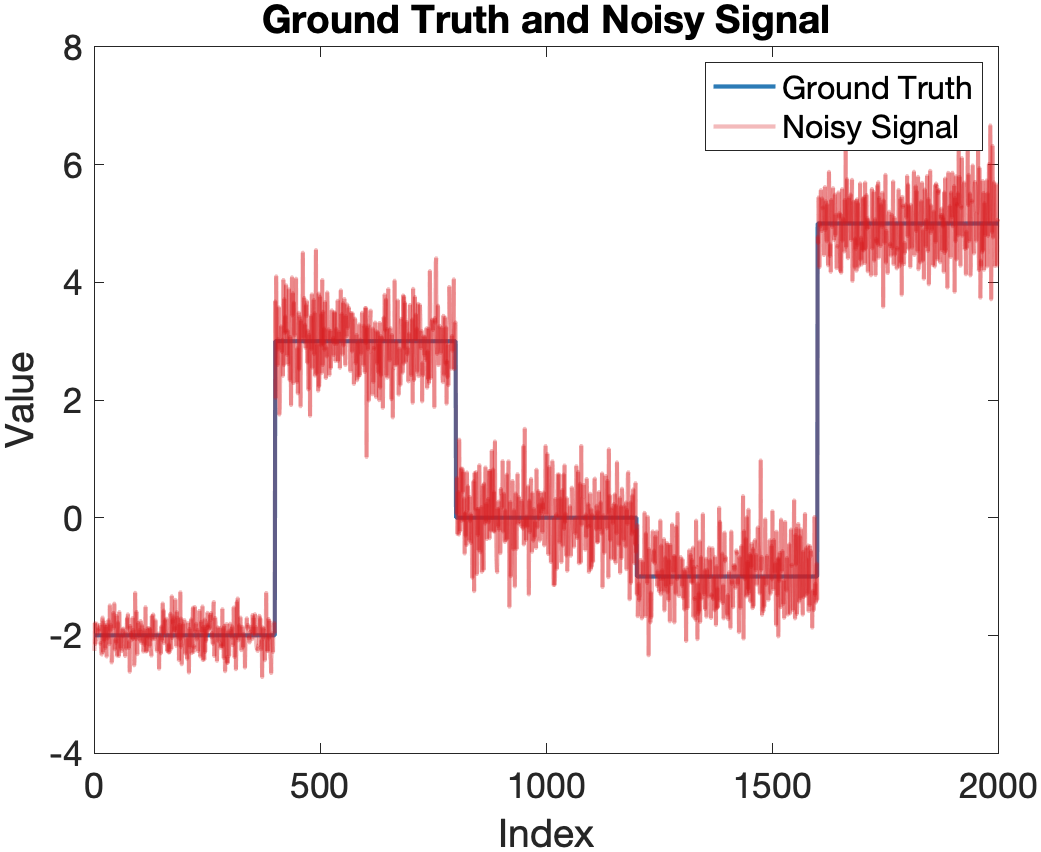}
\caption{The ground truth and noisy signals used for the one-dimensional total variation denoising problem.}
\label{fig:tv1d_truth}
\end{figure}

\begin{figure}[t]
\centering
\includegraphics[width=0.45\linewidth]{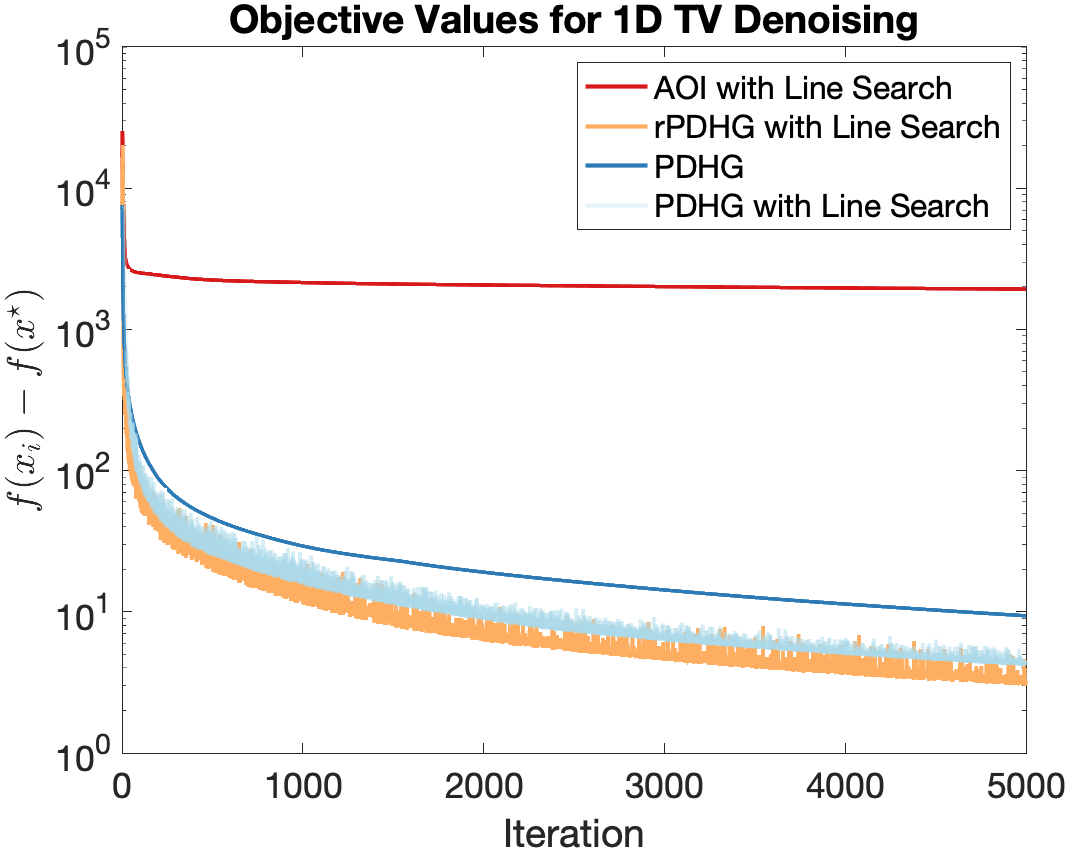}\,\includegraphics[width=0.45\linewidth]{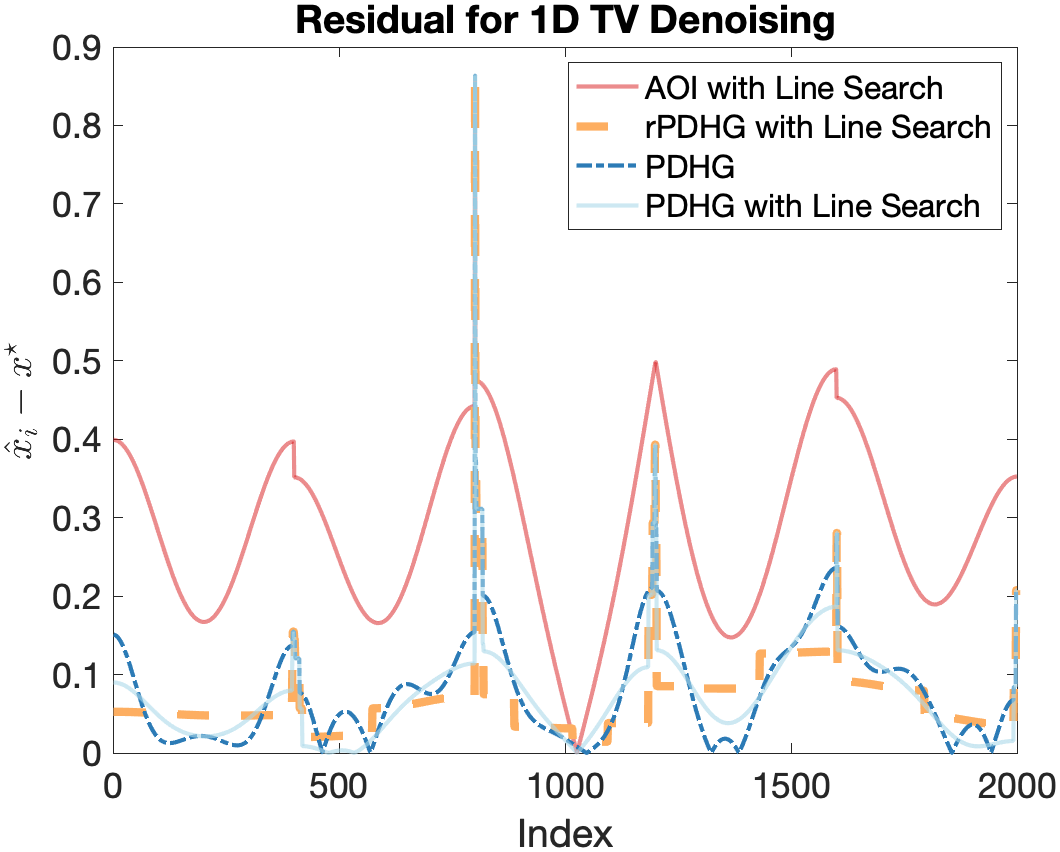}
\caption{(Left) A comparison of the objective values of different algorithms solving the one-dimensional total variation denoising problem. The line search proposed in this manuscript converges with fewer iterations than either constituent line search alone. The objective value jumps often as different step sizes are chosen. The raw objective value for each step is shown for academic purposes; in practice, one would choose the best objective value so far. The PDHG line performs better than standard PDHG. The AOI line search alone converges to a poor solution. (Right) A comparison of the best answers generated by the different optimization algorithms after 1000 iterations, plotted as the residual. }
\label{fig:tv1d_objs}
\end{figure}

Consider the one-dimensional total variation denoising problem (cf. \cite{chambolle2004algorithm}). We take a signal $b$ of length $n$ with additive white noise and solve the problem 
\begin{equation*}
\underset{x \in \mathbb{R}^n}{\text{minimize}}\, \frac{1}{2}\lVert x - b \rVert_2^2 + \lambda\lVert\nabla x \rVert_1,
\label{eq:tv_1d}
\end{equation*}
where $\lambda > 0$ is a chosen parameter and $\nabla: \mathbb{R}^n \to \mathbb{R}^n$ is the one-dimensional finite difference operator: 
\begin{equation*}
\nabla x_i = \begin{cases}\abs{x_1 - x_n}, &i = 1 \\ \abs{x_i - x_{i-1}}, \qquad &i = 2, \ldots, n.  \end{cases}
\end{equation*}
The signal is created by specifying a number of segments $num_{\text{segs}}$ and length of segments $len_{\text{segs}}$ and letting $n =  num_{\text{segs}} * len_{\text{segs}}$. Then $b$ is created by taking $num_{\text{segs}}$ random integers $z \in [-5, 5]$ and creating a step function from $len_{\text{segs}}$. White noise is then added to this to create the measurements $b$. For this experiment, we chose $\lambda = 1$. 

Here $f(x) = \frac{1}{2}\lVert x - b \rVert_2^2$ and $g(x) = \lambda\lVert x \rVert_1$. The $A$ operator for the algorithm is the finite difference operator $\nabla$. Creating $B$ for the line search requires first creating a matrix representation of $A$ and then creating $B$ by choosing $\theta > \frac{1}{\norm{A}^2}$ and solving for $B$ from
\begin{equation*}
BB^* = \frac{1}{\theta} \mathcal{I} - AA^*,
\end{equation*}
via Cholesky factorization. Figure \ref{fig:tv1d_truth} shows the ground truth and the noisy function. Figure \ref{fig:tv1d_objs} shows the objective values for each iteration of the four different algorithms along with the residual of the final answer produced by the algorithms. Note that for the three algorithms that are not rPDHG, we did extensive parameter searches to find the best parameters for apt comparison. For rPDHG, however, we simply set $\beta_0 = 1$ and choose a $\tau_0 > 0$, avoiding a costly parameter search.

\subsection{2D-TV Denoising}\label{sec:2dtv}
\begin{figure}[t]
\centering
\includegraphics[width=0.35\linewidth]{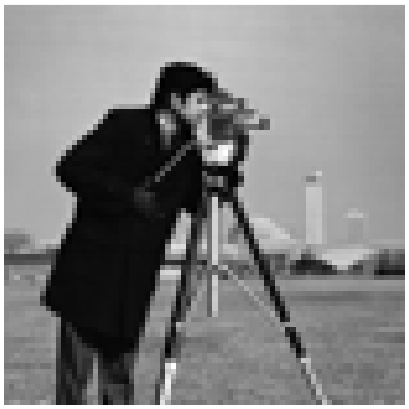}\includegraphics[width=0.35\linewidth]{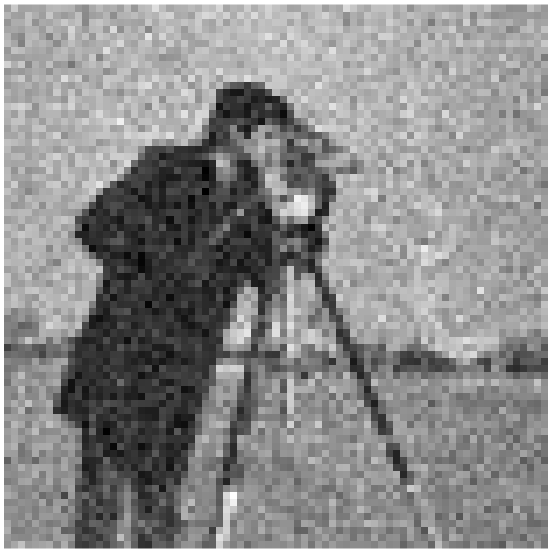}
\caption{The cameraman image used for the two-dimensional total variation denoising problem. On the right is the image with white noise added.}
\label{fig:cameraman}
\end{figure}

The two-dimensional total variation denoising problem follows the standard ROF model \cite{rudin1992nonlinear, chambolle2004algorithm}:
\begin{equation}
\underset{x \in \mathbb{R}^n}{\text{minimize}}\, \frac{1}{2}\lVert x - b\rVert_2^2 + \lVert \nabla x \rVert_{2, 1}.
\label{prob:2dtv}
\end{equation}
The experiments conducted for this manuscript used the \texttt{cameraman.tif} image in MATLAB. The image was read in as a $256\times256$ double and scaled so that its maximum value is $1$. The image was then downsized to $77\times77$ pixels for memory considerations.  (As discussed previously in Section \ref{sec:drawbacks}, for image processing tasks the $B$ matrix can have a large memory footprint. Therefore, to run this example on a desktop computer the image was scaled down.)  Noise was generated as a $77\times77$ matrix drawn from a standard normal distribution and multiplied by $0.08$. The original image and the noised image are shown in Fig. \ref{fig:cameraman}.

The discrete gradient operator $\nabla$ is defined $\nabla: \mathbb{R}^{n \times m} \to \mathbb{R}^{n\times m\times2}$:
\begin{align}
\nabla{x}^1_{i,j} &= 
\begin{cases} 
x_{i+1, j} - x_{i, j}, \qquad &i = 1, \ldots, m-1 \\
x_{1, j} - x_{n, j}, \qquad &i = m
\end{cases} \\
\nabla{x}^2_{i,j} &= 
\begin{cases} 
x_{i, j+1} - x_{i, j}, \qquad &i = 1, \ldots, n-1 \\
x_{i, n} - x_{i, 1}, \qquad &i = n.
\end{cases}
\label{eq:tv_op2}
\end{align}
The $\ell_{2,1}$ norm is defined as
\begin{equation*}
\norm{\nabla x}_{2, 1} = \sum_{i=1}^{m}\sum_{j=1}^n \abs{ (\nabla{x}^1_{i, j})^2  + (\nabla{x}^2_{i, j})^2}^{1/2}.
\end{equation*}

\begin{figure}[t]
\centering
\includegraphics[width=0.5\linewidth]{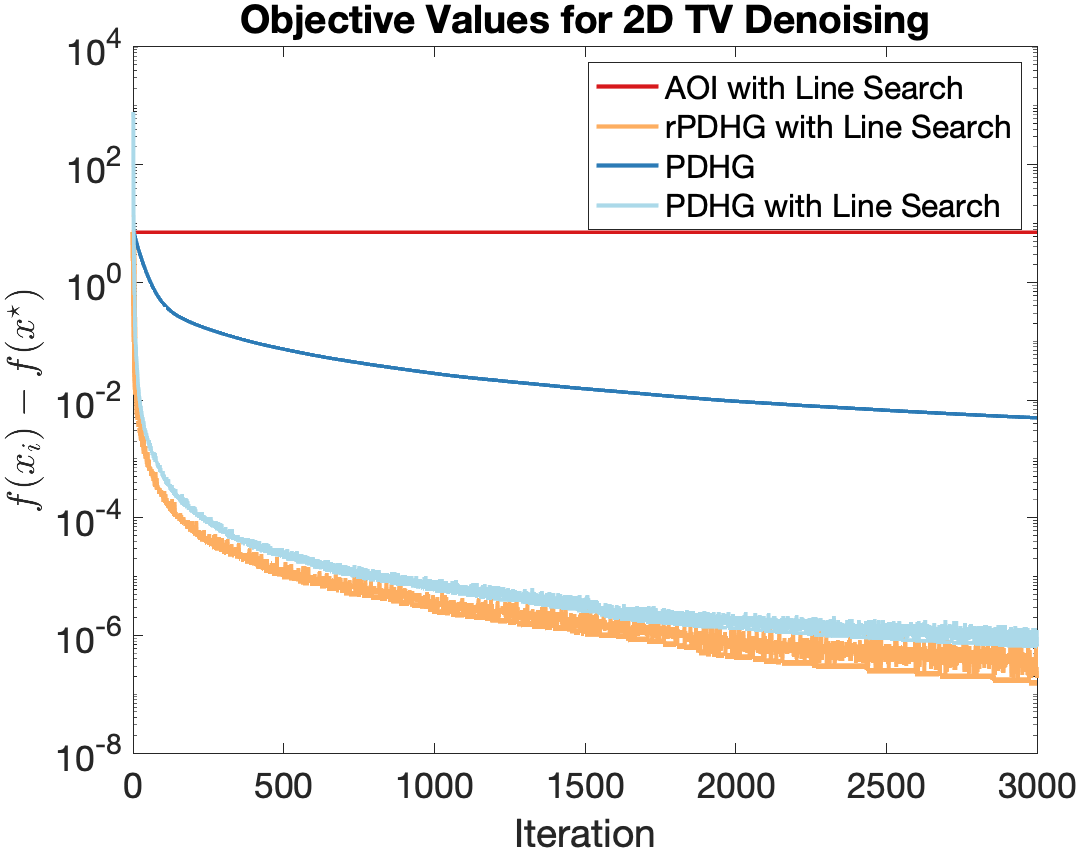}\includegraphics[width=0.4\linewidth]{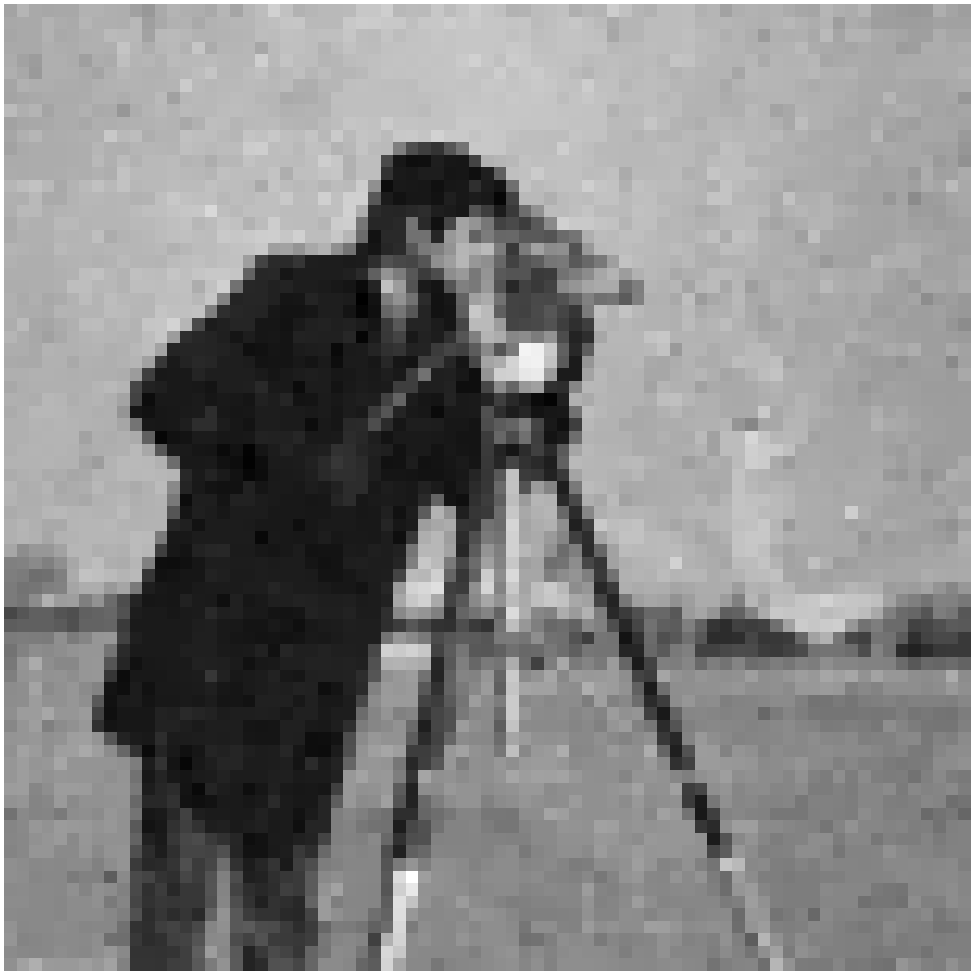}
\caption{(Left) A comparison of the objective values of the different optimization algorithms used for the two-dimensional total variation denoising problem. The AOI line search converged to a poor solution. The line search proposed in this manuscript took fewer iterations to converge to any given level of error and did not require searching over any parameters. The comparison against PDHG and PDHG with line search is done against their best sets of parameters. (Right) The final output of rPDHG, the line search proposed in this manuscript, solving the two dimensional TV denoising problem as in \eqref{prob:2dtv}.}
\label{fig:tv2d_objs}
\end{figure}

Similar to the one-dimensional problem, the discrete gradient operator is our $A$ operator for the minimization problem. To calculate $B$, we first turn $A$ into a matrix from its matrix-free form and then choose $\theta > \frac{1}{\norm{A}^2}$ and solve for $B$ from
\begin{equation*}
BB^* = \frac{1}{\theta} \mathcal{I} - AA^*
\end{equation*}
via Cholesky factorization. A comparison of the objective values for the different algorithms used is in Fig. \ref{fig:tv2d_objs} (left). PDHG and the line search for the AOI formulation of PDDR converge slowly. The line search proposed in this manuscript and the Malitsky line search both tend to ``oscillate'' around a solution. For a fixed number of iterations, our line search and the Malitsky line search have the lowest error. The final output of the proposed algorithm is Fig. \ref{fig:tv2d_objs} (right). 

% \begin{figure}[t]
% \centering
% \includegraphics[width=0.4\linewidth]{figs/2dtv_gpdhg_final.png}
% \caption{The final output of gPDHG, the line search proposed in this manuscript, solving the two dimensional TV denoising problem as in \eqref{prob:2dtv}. The objective value plot is in Fig. \ref{fig:tv2d_objs}.}
% \label{fig:tv2d_final}
% \end{figure}

\section{Application to Magnetic Resonance Imaging Reconstruction}
\label{sec:mri}
\begin{figure}[t]
\begin{centering}
\includegraphics[width=0.8\linewidth]{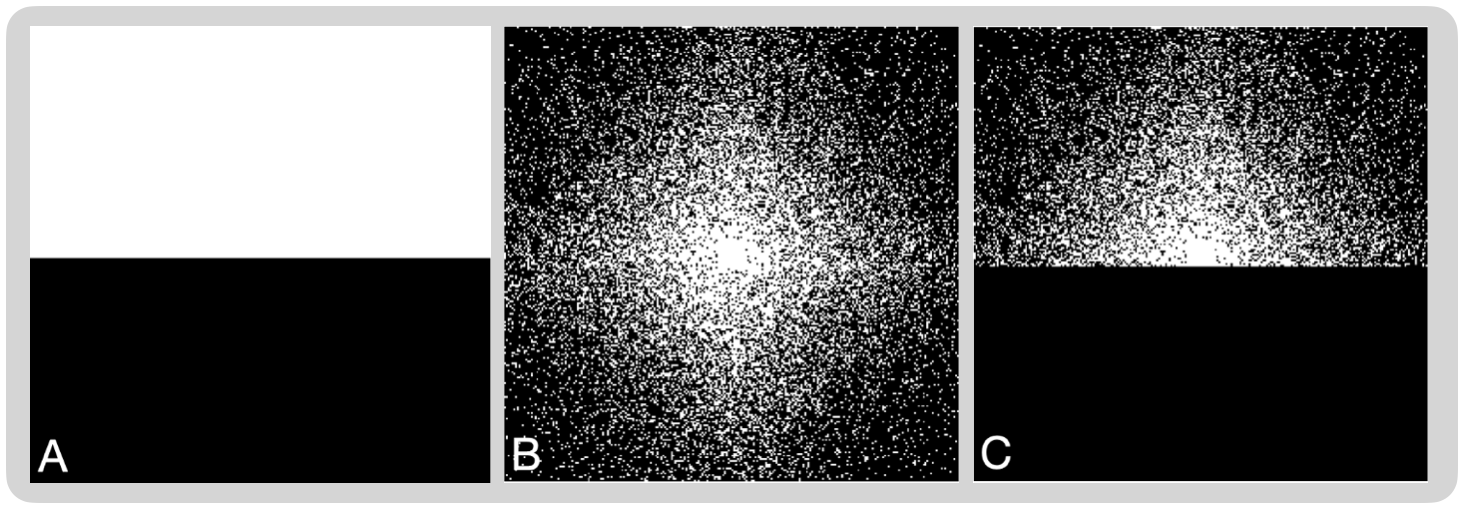}
\caption{ Sampling masks for three different methods of acceleration, where white indicates collected samples. A) is the sampling mask associated with partial Fourier sampling, where the entire top half plus some number of lines around the origin are collected. B) is the sampling mask associated with compressed sensing using a variable density Poisson disc sampling mask \cite{dwork2021fast}. C) is the sampling mask for the novel combination of partial Fourier and compressed sensing reconstruction presented in this manuscript. It is the intersection of masks A) and B).}
\label{fig:sampling_mask_cspf}

\end{centering}
\end{figure}

The driving force behind this line search is an application to accelerated magnetic resonance imaging (MRI) reconstruction. We present a novel combination of acceleration methods, compressed sensing and partial Fourier with homodyne detection. This problem can be written in the form of \eqref{eq:prob_def}.

\subsection{MRI Preliminaries}
An MRI machine is a Fourier sensing device. A simple acquisition model the data from the machine $b$ as
\begin{equation}\label{mri:acq1}
    b = \mathcal{F}\left\{x\right\} + \text{noise} \,
\end{equation}
where $\mathcal{F}$ represents the 2 dimensional Fourier transform, $x$ is the subject being imaged, and there is noise in the acquisition process \cite{Nishimura2010}. If a sufficiently large amount of data is collected to satisfy the Nyquist-Shannon sampling theorem independently in each direction, the the image can be reconstructed with an inverse Discrete Fourier Transform: $\hat{x}=F^{-1}\,b$.

The scan time of MRI is directly proportional to the number of samples acquired. Acquiring a fully sampled grid takes a long time; to accelerate the acquisition, one can reduce the number of samples. That is, there is some data mask $D$ (a matrix comprised of a subset of rows of the identity matrix) such that 
\begin{equation*}\label{mri:acq2}
    b = D\,Fx + \text{noise}\in \mathbb{C}^{nmp},
\end{equation*}
where the data is acquired on an $n \times m$ Cartesian grid and $p \in (0, 1]$ represents the undersampling factor corresponding to the number of $1$s on the diagonal of $D$.

Reconstructing the image given the acquired samples is then a noisy, underdetermined inverse problem. Common approaches to solving this problem include compressed sensing \cite{lustig2007sparse, candes2008introduction, dwork2021utilizing} and machine learning techniques \cite{sandino2020compressed, knoll2020deep}.

Matrix-free implementations of the algorithm and problem statement are desirable due to memory limitations and to take advantage of the Fast Fourier Transform. A parameter-free version of PDHG is also desirable so the MRI technologist does not have to set a parameter by hand which can prevent the reconstruction of a high-quality image within a reasonable number of iterations. 

\subsection{Compressed Sensing with Homodyne Detection}

Compressed sensing incorporates the \textit{a priori} knowledge that natural images are approximately sparse after an appropriate linear transformation (e.g., a discrete Wavelet transform ) to reconstruct high-quality images with from fewer samples than those that would be required by the Nyquist-Shannon sampling limit \cite{lustig2008compressed, candes2008introduction}.  Partial Fourier reconstruction incorporates the fact that the Fourier transform of a real image is Hermitian symmetric into the reconstruction process \cite{mcgibney1993quantitative, bydder2005partial}. Though MR images are not generally real due to variances in the $B_0$ magnetic field and non-linearities in the $B_1$ magnetic field, the phase is almost all low-bandwidth.  Partial Fourier with homodyne detection takes advantage of this fact \cite{noll1991homodyne} by estimating the phase from a low-pass filtered version of the image. Once estimated, phase correction makes the data real and the assumptions of partial Fourier reconstruction are satisfied.

\begin{figure}[t]
\begin{centering}
\includegraphics[width=0.95\linewidth]{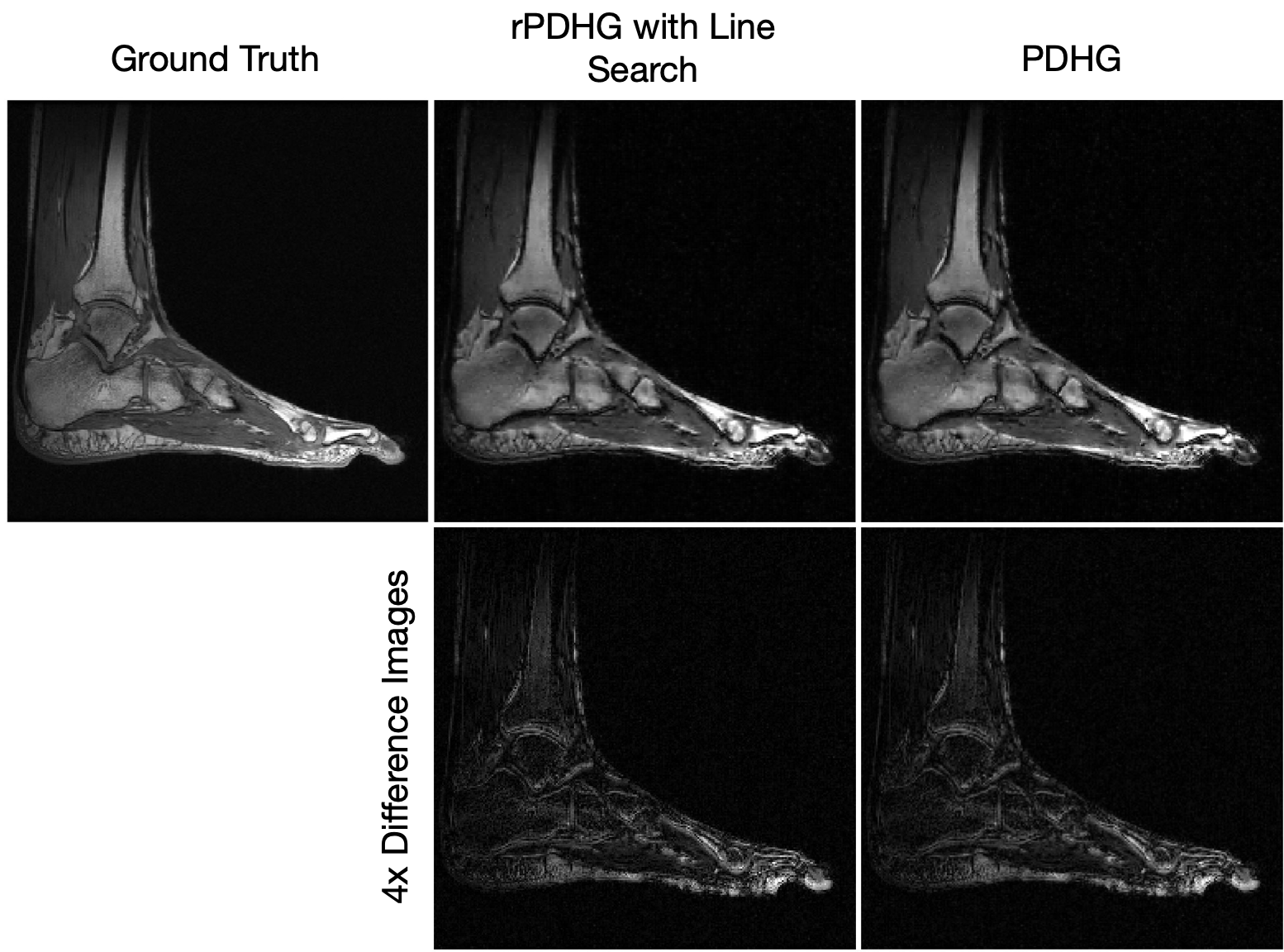}
\caption{Comparison of a fully-sampled reconstruction of a sagittal slice of an ankle to two methods of reconstruction from data undersampled according to sampling mask Fig. \ref{fig:sampling_mask_cspf}-C. (Center) Reconstruction of the undersampled ankle data using the formulation of the problem in \eqref{eq:ch4_cspf_opt_problem_ind}, solved with the line search presented in this manuscript. (Right) This reconstruction was made using the formulation of the problem in \eqref{eq:ch4_cspf_opt_problem_ind}, solved with the PDHG method. Bottom row shows 4 times the magnitude difference between the accelerated reconstructions and the fully-sampled reconstruction.}
\label{fig:ankle_all}
\end{centering}
\end{figure}

\begin{figure}[!t]
\centering
\includegraphics[width=0.95\linewidth]{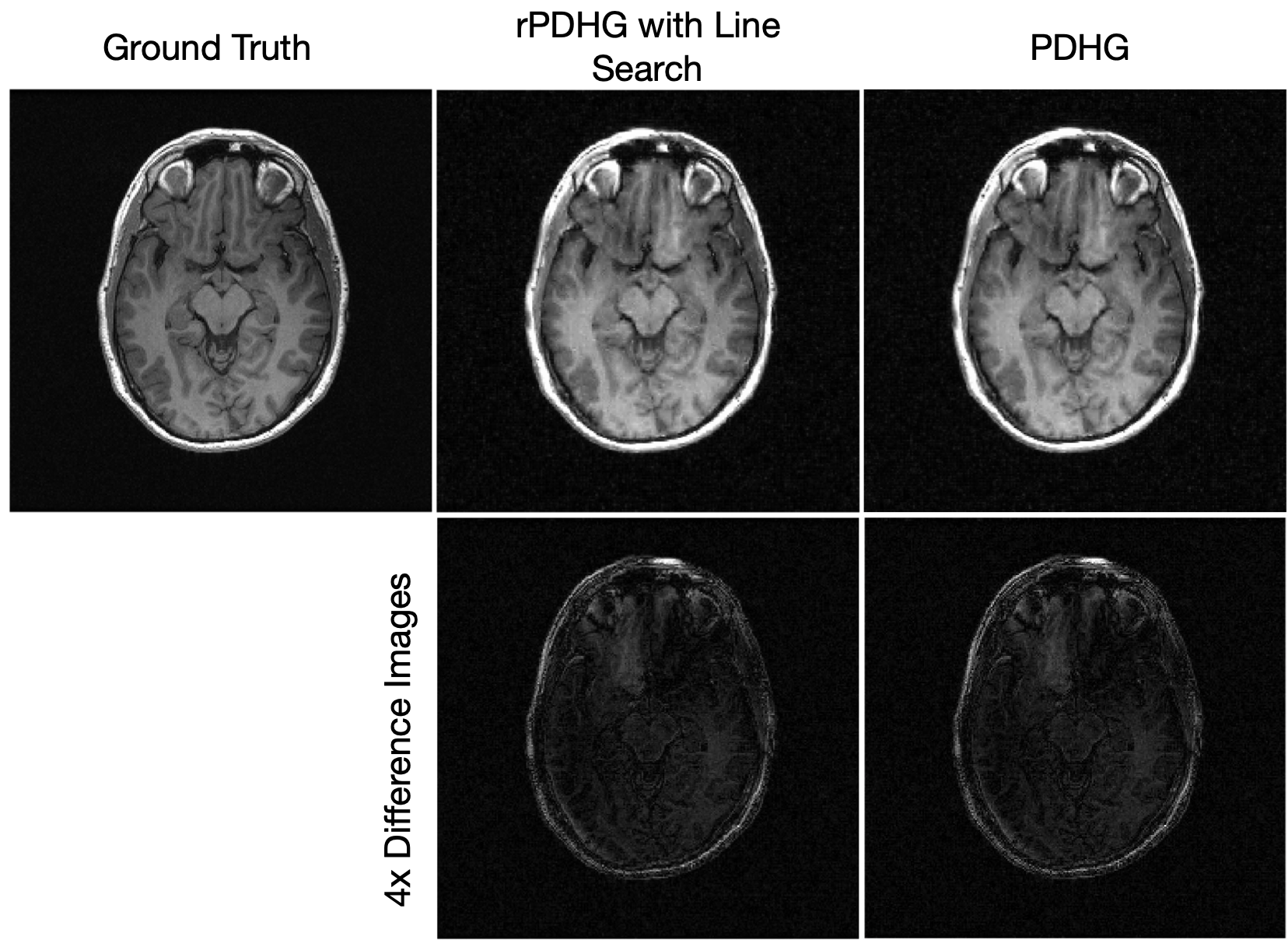}
\caption{Comparison of a fully-sampled reconstruction of an axial slice of a brain to reconstruction from data undersampled according to sampling mask Fig. \ref{fig:sampling_mask_cspf}-C. (Center) Reconstruction of the undersampled ankle data using the formulation of the problem in \eqref{eq:ch4_cspf_opt_problem_ind}, solved with the line search presented in this manuscript. (Right) Reconstruction of the undersampled ankle data using the formulation of the problem in \eqref{eq:ch4_cspf_opt_problem_ind}, solved with the standard PDHG method given in \eqref{eq:iterations-pdhg}. Bottom row shows 4 times the magnitude difference between the accelerated reconstructions and the fully-sampled reconstruction. The reconstruction using the new line search shows the same quality reconstruction without requiring a lengthy parameter search.}
\label{fig:brain_cs_ls}
\end{figure}

% \begin{figure}[t]
% \begin{centering}
% \includegraphics[width=0.75\linewidth]{figs/ankle_compare.png}
% \caption{The two undersampled reconstructions from Fig. \ref{fig:ankle_all} zoomed in on the toe. There is increased detail around the bones and edges of the subject in the reconstruction using gPDHG with linesearch. The reconstruction using PDHG exhibits blurring near the edge of the subject and around the bone structure. }
% \label{fig:ankle_recon_zoomed}
% \end{centering}
% \end{figure}

Partial Fourier sampling collects data with a sampling pattern consisting of the top $\nu$ portion, where $\nu>0.5$. Common values of $\nu$ are $\nu = 5/8$ and $\nu = 9/16$. Homodyne detection uses a low-pass filtered version of the image to estimate the phase, $\Phi$. Once the phase is estimated, the image is reconstructed with homodyne detection using
\begin{equation*}
  x = P_\Phi\,b = Re\left[ \Phi \, F^{-1} \, R \, b \right],
  \label{eq:compressedSensing}
\end{equation*}
where $R$ is a weighting ramp that smoothly transitions from $2$ to $0$, $\Phi=\exp\left( -i\,\text{angle}(M_s) \right)$, $M_s=F^{-1}\,L\,b$, and $L$ is a low-pass rectangular filter \cite{pauly2005partial}.

The crucial realization is that once the phase $\Phi$ is estimated, the $P_{\Phi}$ operator---which implements reconstruction with homodyne detection---is a linear transformation with a computationally efficient adjoint. This allows us to combine the phase estimate $\Phi$ into a compressed sensing image reconstruction as the following optimization problem:
\begin{equation*}
  \underset{\xi}{\text{minimize}} \hspace{0.5em} \| \Psi\,P_{\Phi}\,\xi \|_1 \hspace{0.5em} \text{subject to} \hspace{0.5em} \|D\,F\,\xi - b\|_2\leq\epsilon,
\end{equation*}
where $\xi$ is a vector of the top $\nu$ Fourier values for the reconstructed image, $\Psi$ is a sparsifying transformation, and $\varepsilon > 0$ is a noise bound. A typical sparsifying transformation is a wavelet transformation. We write this as an unconstrained problem by using an indicator function to enforce data consistency:
\begin{equation}
  \xi^\star \in \underset{\xi}{\text{argmin}} \hspace{0.5em} \mathbb{I}\left(M\,\xi\in B_{\varepsilon}[b]\right) + \lVert \Psi P_{\Phi}\,\xi\rVert_1,
  \label{eq:ch4_cspf_opt_problem_ind}
\end{equation}
where $\mathbb{I}$ is an indicator function that requires the Fourier values $\xi$ at the known, collected indices specified by sampling mask $D$ to be within $\varepsilon$ of the collected data $b$. Setting $\varepsilon = 0$ enforces data consistency completely. The second term contains the sparsifying transform $\Psi$ and the homodyne operator $P_\Phi$. By letting $A=\Psi P_{\Phi}$, $f(\xi)=\mathbb{I}\left(M\,\xi\in B_{\varepsilon}[b]\right)$ and $g(\xi)=\lVert \Psi P_{\Phi}\,\xi\rVert_1$, we see that this problem is of the form \eqref{eq:prob_def} so we solve it via rPDHG with line search. Once $f^\star$ is determined by solving \eqref{eq:ch4_cspf_opt_problem_ind}, the image is reconstructed according to $x^\star = P_\Phi \xi^\star$.

\subsubsection{Matrix-free implementation}
To aid in analysis, it will help to write $A = \Psi P_{\Phi}$ in its full, matrix-form implementation:
% \begin{equation}
%     \label{eq:A_matrix}
%     A = \Psi \begin{bmatrix} \mathcal{I} & 0  \end{bmatrix}\begin{bmatrix}\Phi_{R} & -\Phi_I \\ \Phi_I & \Phi_R \end{bmatrix}\begin{bmatrix}F_R & F_I \\ -F_I & F_R \end{bmatrix}\begin{bmatrix}R&\\&R\end{bmatrix}\begin{bmatrix}D^*&\\&D^* \end{bmatrix}.
% \end{equation}
\begin{equation}
    \label{eq:A_matrix_2}
    A = \Psi \circ \text{Re} \circ \Phi \circ F^{-1} \circ R \circ D^*.
\end{equation}
Here $D^* : \mathbb{C}^{nmp} \to \mathbb{C}^{n\times m}$ is a zero-filling operator. The linear transformation $D$ selects only the indices at which data was collected (non-zero in the sampling mask) so $D^*$ returns a zero-filled full-sized image. The linear transformation $R : \mathcal{C}^{n\times m} \to \mathcal{C}^{n\times m}$ applies a weighting matrix, $\mathcal{F} : \mathbb{C}^{n\times m} \to \mathbb{C}^{n\times m}$ is the two dimensional discrete Fourier transform, $\Phi : \mathbb{C}^{n\times m} \to \mathbb{C}^{n\times m}$ applies the phase estimate. $\text{Re} : \mathbb{C}^{n\times m} \to \mathbb{R}^{n \times m}$ takes the real part, and $\Psi : \mathbb{R}^{n\times m} \to \mathbb{R}^{n\times m}$ is the orthogonal discrete wavelet transform.

% All of the matrices in the above equation should be understood as diagonal matrices, as $A: \mathbb{R}^{2np} \to \mathbb{R}^{2nm}$ operates on the undersampled data vector $b$. The $D^*$ operator is a zero-filling operator. $D$ selected only the indices at which data was collected (non-zero in the sampling mask) so $D^*$ returns to a full-sized image, represented as a diagonal matrix. $R$ is a diagonal matrix which applies a weighting ramp for the homodyne reconstruction. $F_R$ and $F_I$ correspond to the real and imaginary parts of the discrete Fourier transform, respectively. $\Phi_R$ and $\Phi_I$ correspond to the real and imaginary parts of the phase of the image. \todo{figure out dimensions}

If $A$ were implemented as a matrix, we could choose $\frac{1}{\theta} \geq \lVert A\rVert^2$ and let $B = \left(\frac{1}{\theta} \mathcal{I} - AA^*\right)^{1/2}$, where $\left(\cdot\right)^{1/2}$ indicates taking the Cholesky factorization of a matrix. When we choose $\frac{1}{\theta}$ that satisfies $\frac{1}{\theta} \geq \lVert A\rVert^2$, then $\frac{1}{\theta} \mathcal{I} - AA^*$ is symmetric and positive definite and the Cholesky factorization is well-defined.  However, for this application, the implementation of $A$ is matrix-free---we write it out in the form of \eqref{eq:A_matrix_2} for analysis.  Note that $AA^* = \Psi P_{\Phi}P_{\Phi}^*\Psi^*$ is \textit{not} equal to a scaled identity. Therefore, to use the line search presented above for an image of reasonable size we need a matrix-free implementation of $B$ such that $AA^* + BB^* = \frac{1}{\theta}\mathcal{I}.$

Let $B$ be defined as
% \begin{equation*}
%     B = \Psi \begin{bmatrix} \mathcal{I} & 0  \end{bmatrix}\begin{bmatrix}\Phi_{R} & -\Phi_I \\ \Phi_I & \Phi_R \end{bmatrix}\begin{bmatrix}F_R & F_I \\ -F_I & F_R \end{bmatrix}\begin{bmatrix}R&\\&R \end{bmatrix}\begin{bmatrix}\left(\frac{1}{\theta} \mathcal{I} - RD^*DR\right)^{1/2}&\\&\left(\frac{1}{\theta} \mathcal{I} - RD^*DR\right)^{1/2} \end{bmatrix}.
% \end{equation*}

\begin{equation}
    \label{eq:B_matrix_2}
    B = \Psi \circ \text{Re} \circ \Phi \circ F^{-1} \circ R \circ Q.
\end{equation}
$Q: \mathbb{C}^{n\times m} \to \mathbb{C}^{n\times m}$ is defined as $\left(\frac{1}{\theta}\mathcal{I} - RD^*DR\right)^{1/2}$. The operator $RD^*DR : \mathbb{C}^{n\times m} \to \mathbb{C}^{n\times m}$ is symmetric so $Q$ is the symmetric square root. In practice, this square root is cheap to compute as the matrix is diagonal.  This form of $B$ was found by writing $A$ in matrix form as in \eqref{eq:A_matrix_2}, forming $AA^*$ and setting $BB^* = \frac{1}{\theta}\mathcal{I} - AA^*$. Note that $R$ is symmetric. Similar to above, we write $B$ out in its matrix form to aid in analysis, but $B$ has a natural matrix-free implementation. With this definition of $B$, and understanding that we are using unitary Fourier and wavelet transformations, we can apply the primal-dual form of Douglas-Rachford to the problem \eqref{eq:DRtoPDHG} and recover the PDHG iterations.

\subsection{Results}
We first compare our proposed line search method against standard PDHG and PDHG with line search \cite{malitsky2018first}. A variable density sampling mask with a sampling burden of $0.3$ is chosen. This corresponds to Fig. \ref{fig:sampling_mask_cspf}-B. Intersecting with the partial Fourier mask as in Fig. \ref{fig:sampling_mask_cspf}-C yields an \textit{effective} sampling burden of approximately $0.17$.  (I.e., the number of samples of the reduced sampling pattern divided by those required to satisfy the Nyquist-Shannon sampling theorem was approximately $0.17$.)  We solve the problem \eqref{eq:ch4_cspf_opt_problem_ind} using each method with the same sampling mask. Each method is run for 4000 iterations and the image is reconstructed using the method of Roemer \cite{roemer1990nmr}. Comparing the reconstructions on ankle data in Fig. \ref{fig:ankle_all} center and right shows that the proposed line search yields the same quality reconstruction without a lengthy parameter search. The difference images show that the primary differences are down in the fine structure near the toes and along other areas where the phase is changing rapidly. 

Fig. \ref{fig:brain_cs_ls} is an example using brain data with the same undersampling pattern applied. On the right is the reconstruction using standard PDHG. In the center is the novel combination of compressed sensing with homodyne detection solve using the line search presented in this paper. Similar to the ankle data, the same quality reconstruction is achieved without a parameter search.

\begin{figure}[t]
\begin{centering}
\includegraphics[width=0.45\linewidth]{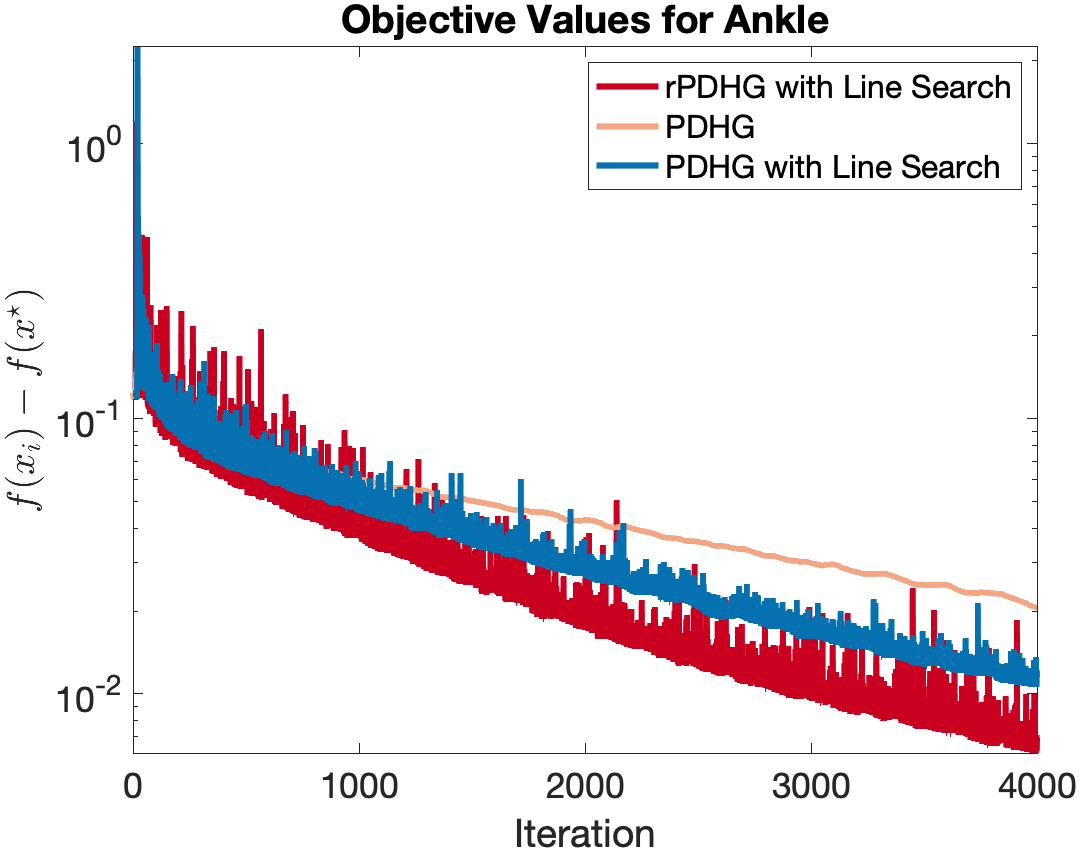}
\includegraphics[width=0.45\linewidth]{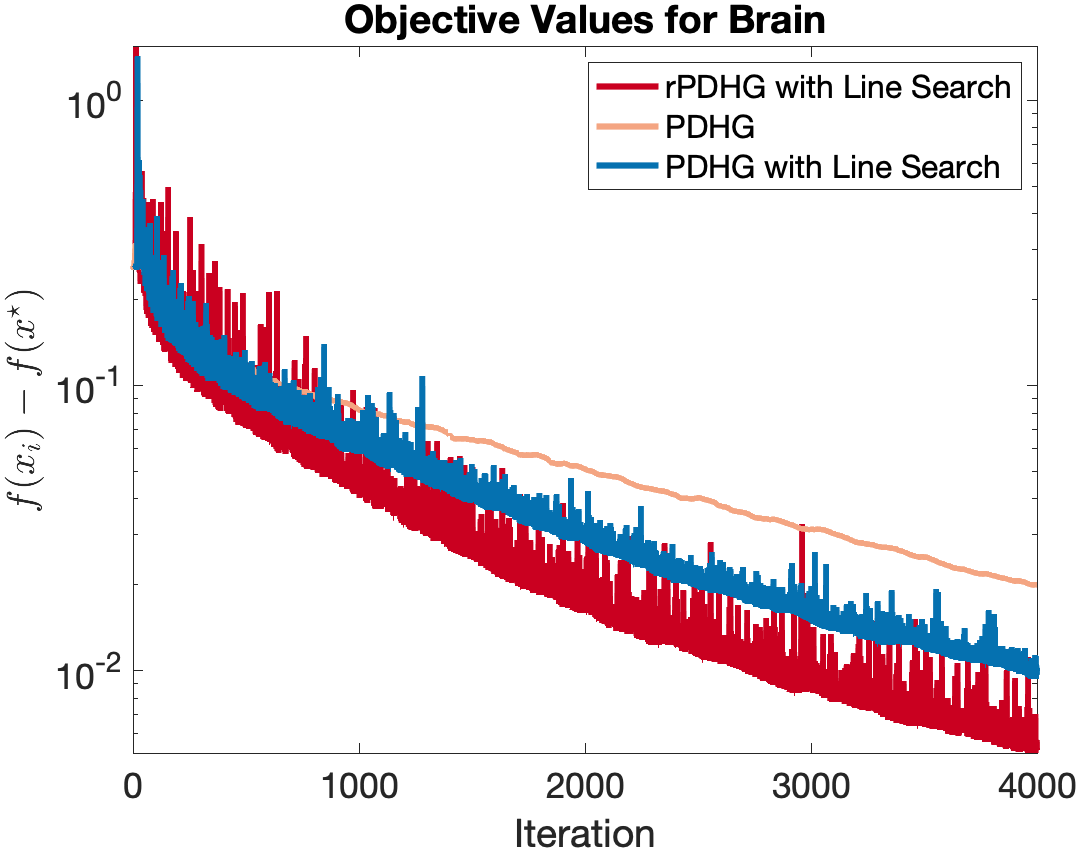}
\caption{Objective value plots for the reconstruction of the same ankle data as Fig. \ref{fig:ankle_all} (left) and the brain data as in Fig. \ref{fig:brain_cs_ls} (right). This compares the performance of rPDHG algorithm with standard PDHG and PDHG with Malitsky's line search. Over time, the line search finds a better solution. The line search jumps in objective value a lot as it chooses different step sizes. The objective value at each iteration is shown for academic purposes, though in practice at each iteration the best so far would be chosen. }
\label{fig:ankle_ls_nols}
\end{centering}
\end{figure}

A comparison of the objective values for reconstruction of the ankle data with and without the line search is presented in Fig. \ref{fig:ankle_ls_nols}. Let $f_o$ be the objective function the algorithm is minimizing (e.g. \eqref{eq:ch4_cspf_opt_problem_ind}). For these plots, let $x^{\star}$ be the optimal point obtained after letting the algorithm run for several thousand iterations. The objective value plotted in Fig. \ref{fig:ankle_ls_nols} is $f_o(x) - f_o(x^{\star})$, or the difference between the current iteration and optimal. As mentioned previously, there is no guarantee that each iterate is feasible. To calculate the objective value we first use the proximal operator of the indicator function for the constraint to make sure the iterate is feasible and then calculate the objective value. The data was run with three different methods: standard PDHG without the line search, PDHG with the Malitsky line search, and relaxed PDHG with the line search and the activation criterion discussed in Sec. \ref{sec:ls_activation}. The objective of the relaxed PDHG algorithm jumps as different step sizes are accepted and converges to a solution better than the PDHG algorithm. The objective values for all iterations are shown for academic purposes; in practice, at any iteration, the result with the best objective value so far would be chosen. The best PDHG run is shown after a lengthy parameter search. The rPDHG run requires no such parameter search. Table \ref{tab:timing} is an example of the timing difference given the parameter search. The time to 1000 iterations is longer for rPDHG but the pretraining time to find the correct set of parameters is zero. 

\begin{table}[!h]
    \centering
    \begin{tabular}{|c||c|c|c|}
    \hline
        Algorithm & Time of 1000 Iterations (s) & Time of Pretraining (s) & Total time (s) \\
        \hline
        PDHG & 38.44 & 653.53 & 691.97 \\
        \hline
        PDHG with Line Search & 98.37 & 1672.29 & 1770.66 \\
        \hline
        AOI with Line Search & 86.67 & 1473.39 & 1560.06 \\
        \hline
        rPDHG with Line Search & 417.8 & 0 & 417.8\\
        \hline
    \end{tabular}
    \caption{An example of the total timing difference between the other algorithms tested and rPDHG with line search given the need for a parameter search for other algorithms. The time to 1000 iterations given a set of parameters is longer for rPDHG with line search, but the total time is less.}
    \label{tab:timing}
\end{table}

\section{Conclusion}

A line search for the primal-dual hybrid gradient method was proposed. Motivated by the connection between the Douglas-Rachford algorithm and the primal-dual hybrid gradient method, this line search expresses PDHG as an averaged operator iteration. An existing line search for averaged operator iterations is then applied to the relaxation parameter. This is combined with an existing line search over the proximal step sizes in PDHG. This new combination of line searches is essentially parameter free, allowing for quality solutions without a lengthy parameter search. Several numerical examples are demonstrated showing the convergence performance. The efficacy of the line search is presented on a novel combination of magnetic resonance imaging reconstruction techniques.

\section*{Acknowledgments}

The authors thank Daniel O'Connor for many useful conversations related to the work of this manuscript.

\newpage
\appendix
\newpage \label{appendix_a}

\section{Computation of B for the MRI Problem}
For the compressed sensing with homodyne reconstruction problem, $A$ is given as
\begin{equation}
    \label{eq:app_a_1}
    A = \Psi \circ \text{Re} \circ \Phi \circ F^{-1} \circ R \circ D^*.
\end{equation}

To use relaxed PDHG with line search, we need a $B$ such that $AA^* + BB^* = \frac{1}{\theta}\mathcal{I}$ for some $\theta > 0$. From \eqref{eq:app_a_1}, we need
\begin{align}
    \left(\Psi \circ \text{Re} \circ \Phi \circ F^{-1} \circ R \circ D^*\right)\circ\left(\Psi \circ \text{Re} \circ \Phi \circ F^{-1} \circ R \circ D^*\right)^* + BB^* &= \frac{1}{\theta}\mathcal{I}\nonumber\\
    \left(\Psi \circ \text{Re} \circ \Phi \circ F^{-1} \circ R \circ D^*\right)\circ\left(D\circ R^* \circ (F^{-1})^* \circ \Phi^* \circ \text{Re}^* \circ \Psi^* \right) + BB^* &= \frac{1}{\theta}\mathcal{I}.
\label{eq:app_a_2}
\end{align}
If we write $B = \Psi\circ\text{Re}\circ\Phi\circ F^{-1}\circ Q$ for unknown $Q$, we can write \eqref{eq:app_a_2} as
\begin{equation*}
\begin{split}
    &\underbrace{\left(\Psi \circ \text{Re} \circ \Phi \circ F^{-1} \circ R \circ D^*\right)}_{A}\circ\underbrace{\left(D\circ R^* \circ (F^{-1})^* \circ \Phi^* \circ \text{Re}^* \circ \Psi^* \right)}_{A^*} \\+ &\underbrace{\left(\Psi \circ \text{Re} \circ \Phi \circ F^{-1} \circ Q\right)}_{B}\circ\underbrace{\left(Q^* \circ (F^{-1})^* \circ \Phi^* \circ \text{Re}^* \circ \Psi^* \right)}_{B^*} = \frac{1}{\theta}\mathcal{I}.
\end{split}
\end{equation*}
Factoring out the similar terms yields
\begin{equation}
\left(\Psi \circ \text{Re} \circ \Phi \circ F^{-1}\right) \circ \left(RD^*DR^* + QQ^*\right)\circ \left((F^{-1})^* \circ \Phi^* \circ \text{Re}^* \circ \Psi^*\right) = \frac{1}{\theta}\mathcal{I}.
\label{eq:app_a_3}
\end{equation}

We are using the unitary Fourier transform $F$ and orthogonal wavelet transform $\Psi$. The phase estimate $\Phi$ is unitary---it is a diagonal matrix of unit complex exponentials. The transpose of the Re operator is $\text{Re}^* : \mathbb{R}^{n\times m} \to \mathbb{C}^{n\times m}$ and maps $a \to a + 0i$, so $\text{Re}\circ\text{Re}^* : \mathbb{R}^{n \times m} \to \mathbb{R}^{n \times m}$ is the identity operator. Therefore if $RD^*DR^* + QQ^* = \frac{1}{\theta}\mathcal{I}$ then \eqref{eq:app_a_3} is satisfied. Noting that $R$ is Hermitian, this condition becomes
\begin{equation}
QQ^* = \frac{1}{\theta}\mathcal{I} - RD^*DR.
\label{eq:app_a_4}
\end{equation}
The matrix $\frac{1}{\theta}\mathcal{I} - RD^*DR$ is diagonal so taking a symmetric square root is computationally efficient.
% \input{app1}
% \input{app3_pddrpdhg}
%\include{app4_ops}
%\include{app2_ls}

%\bibliographystyle{unsrt}
%\bibliography{references}  %%% Uncomment this line and comment out the ``thebibliography'' section below to use the external .bib file (using bibtex) .

\end{document}